\documentclass[10pt]{article}
\usepackage{amssymb,amsmath,amsfonts,amsthm,enumerate,color}
\usepackage[toc,page,title,titletoc,header]{appendix}
\usepackage{graphicx}
\usepackage{indentfirst}
\usepackage{multicol}
\usepackage{booktabs}
\usepackage{url}

\numberwithin{equation}{section}

\newtheorem{theorem}{Theorem}

\newtheorem{lemma}{Lemma}

\def \Vh0{\stackrel{\circ}{V}_h} \def\to{\rightarrow}

\newcommand{\lc}
{\mathrel{\raise2pt\hbox{${\mathop<\limits_{\raise1pt\hbox
{\mbox{$\sim$}}}}$}}}

\newcommand{\gc}
{\mathrel{\raise2pt\hbox{${\mathop>\limits_{\raise1pt\hbox{\mbox{$\sim$}}}}$}}}

\newcommand{\ec}
{\mathrel{\raise2pt\hbox{${\mathop=\limits_{\raise1pt\hbox{\mbox{$\sim$}}}}$}}}

\def\bb{\begin{equation}} \def\ee{\end{equation}}

\def\beqn{\begin{eqnarray}}  \def\eqn{\end{eqnarray}}

\def\beqnx{\begin{eqnarray*}} \def\eqnx{\end{eqnarray*}}

\def\bn{\begin{enumerate}} \def\en{\end{enumerate}}

\def\bd{\begin{description}} \def\ed{\end{description}}
\DeclareMathOperator*{\argmin}{arg\,min}


\newcommand{\EE}{\mathbb{E}}

\begin{document}

\title{Asynchronous Coordinate Descent under More Realistic Assumptions }
\author{Tao Sun\thanks{
College of Science, National University of Defense Technology,
Changsha, 410073, Hunan,  China. Email: \texttt{nudttaosun@gmail.com;nudtsuntao@163.com} }
\and Robert Hannah\thanks{
  Department of Mathematics, UCLA, 601 Westwood Plz, Los Angeles, CA 90095, USA, Email: \texttt{RobertHannah89@math.ucla.edu}}
  \and   Wotao Yin\thanks{
  Department of Mathematics, UCLA, 601 Westwood Plz, Los Angeles, CA 90095, USA, Email: \texttt{wotaoyin@math.ucla.edu}}
}

\maketitle

\begin{abstract}
Asynchronous-parallel algorithms have the potential to vastly speed up algorithms by eliminating costly synchronization. However, our understanding of these algorithms is limited because the current convergence of asynchronous (block) coordinate descent algorithms are based on somewhat unrealistic assumptions. In particular, the age of the shared optimization variables being used to update a block is assumed to be independent of the block being updated. Also, it is assumed that the updates are applied to randomly chosen blocks. In this paper, we argue that these assumptions either fail to hold or will imply less efficient implementations.
%In this paper, we first briefly present experiments that show that they are often highly dependent.

We then prove the convergence of asynchronous-parallel block coordinate descent under more realistic  assumptions, in particular, always without the independence assumption. The analysis permits both the deterministic (essentially) cyclic and random rules for block choices. Because a bound on the asynchronous delays may or may not be available, we establish convergence for both bounded delays and unbounded delays. The analysis also covers nonconvex, weakly convex, and strongly convex functions. We construct Lyapunov functions that directly model both objective progress and delays, so delays are not treated errors or noise. A continuous-time ODE is provided to explain the construction at a high level.
\end{abstract}

%\textbf{Keywords: Alternating minimization algorithms, Bregman distance, nonconvex, Kurdyka-{\L}ojasiewicz property}

\section{Introduction}\label{intro}
%%%%%%%%%%%%%%%%%%%%%%%%%%%%%%%%%%%%%%%%
In this paper, we consider the asynchronous-parallel block coordinate descent (async-BCD) algorithm for solving
\small
\begin{align}
\min_{x\in \mathbb{R}^N}\, f(x) = f(x_1,\ldots,x_N),
\end{align}
\normalsize
where $f$ is a differentiable function whose gradient is $L$-Lipschitz continuous.
%\footnote{It is possible to extend our results to composite minimization (minimizing $f+g$, where $f$ has $L$-Lipschitz gradient, and $g$ is convex and proximable). But, we only consider the simpler case where $g\equiv 0$.}.

Async-BCD \cite{LiuWrightReBittorfSridhar2015_asynchronous,LiuWright2015_asynchronous,peng2016arock} has virtually the same implementation as regular BCD. The difference is that the threads doing the parallel computation will not wait for all others to finish and share their updates, but merely continue to update with the most recent updates available\footnote{Also the step size needs to be modified to ensure convergence results hold. However in practice traditional step sizes appear to work, barring extreme circumstances.}. In traditional algorithms, latency, bandwidth limits, and unexpected drains on resources, that delay the update of even a single thread will cause the entire system to wait. By eliminating this costly idle time, asynchronous algorithms can be much faster than traditional ones. % (see section \ref{ssec:Motivation} for results. We observed a 7x speedup in one scenario, using the same computing resources.).

In async-BCD, each agent continually updates the solution vector, one block at a time, leaving all other blocks unchanged. Each block update is a read-compute-update cycle. It begins with an agent reading $x$ from shared memory or a parameter server and saving it in a local cache as $\hat{x}$. Then, the agent computes a block partial gradient $-\frac{\gamma_k}{L}\nabla_{i} f(\hat{x})$, where $\gamma_k$ is a step size. The computing can start before the reading is completed. If $\nabla_{i} f(\hat{x})$ does not require all components of $\hat{x}$, only the required ones are read. The final step of the cycle depends on the parallel system setup. In a shared memory setup, the agent reads $x_{i}$ again and writes $x_{i}-\frac{\gamma_k}{L}\nabla_{i} f(\hat{x})$ to $x_{i}$. In the server-worker setup, the (worker) agent can send $-\frac{\gamma_k}{L}\nabla_{i} f(\hat{x})$ or just $\nabla_{i} f(\hat{x})$ to the server and let the server update $x_i$. Other setups are possible, too. The iteration counter $k$ increments upon the completion of any block update, and the updating block is denote as $i_k$.

Because the block updates are asynchronous, when a block update is completed, the $\hat{x}$ that is read and then used to compute this update can be outdated at the completion time.

The iteration of asyn-BCD is, therefore, modeled \cite{LiuWrightReBittorfSridhar2015_asynchronous} as
\small
\begin{align}
  x^{k+1}_{i_k}=x^k_{i_k}-\frac{\gamma_k}{L}\nabla_{i_k} f(\hat{x}^k)\label{eq:ARock-algorithm},
\end{align}
\normalsize
%
%where $i_k$ is the index of the coordinate block updated, $\widehat{x}^k$ is the potentially outdated version of the solution vector being used to compute the update, $\gamma_k$ is the step size, and $x^{k}_{i_k}$ and $x^{k+1}_{i_k}$ represent the value of $x_{i_k}$ in the shared memory right before and after the update, respectively.
%In this paper we prove convergence results under a variety of assumptions on the form of $i_k$ and $\widehat{x}^k$, which model a variety of implementation strategies. These will be discussed in section \ref{ssec:Setup}.
and $x_j^{k+1}=x_j^k$ for all non-updating blocks $j\neq i_k$. The convergence behavior of this algorithm depends on the sequence of updated blocks $i_k$, the step size sequence $\gamma_k$, as well as the ages of $\hat{x}^k$ relative to $x^k$, which is also called \emph{delays}. We define the \emph{delay vector} $$\vec{j}(k)=(j(k,1),j(k,2),\ldots,j(k,N))\in \mathbb{Z}^N.$$ %The components of this vector represent how many iterations out-of-date the vector $\widehat{x}^k$ is\footnote{We follow the notation and definitions from \cite{hannah2016unbounded}.}.
More precisely,
\small
\begin{align}\label{def:jk}
    \hat{x}^k=(x^{k-j(k,1)}_1,x^{k-j(k,2)}_2,\ldots,x^{k-j(k,N)}_N).
\end{align}
\normalsize
The \emph{$k$'th delay} (or \emph{current delay}) is $j(k)=\max_{1\leq i\leq N}\{j(k,i)\}$.

%%%%%%%%%%%%%%%%%%%%%%%%%%%%%%%%%%%%%%%%
\subsection{Dependence between delays and blocks}\label{ssec:Motivation}
%%%%%%%%%%%%%%%%%%%%%%%%%%%%%%%%%%%%%%%%
In previous analyses \cite{LiuWright2015_asynchronous,LiuWrightReBittorfSridhar2015_asynchronous,peng2016arock,hannah2016unbounded}, it is assumed that the block index $i_k$ and the delay $\vec{j}(k)$ were independent sequences. This simplifies proofs, for example, giving $\EE_{i_k} ( P_{i_k}\nabla f(\hat{x}^k) ) = \tfrac{1}{N} \nabla f(\hat{x}^k)$ when $i_k$ is chosen at random, where $P_i$ denotes the projection to the $i$th block. Without independence, $\vec{j}(k)$ will depend on $i_k$, causing $\hat{x}^k$ to be different for each possible $i_k$ and breaking the equality. The independence assumption is unrealistic in practice. Consider a problem where some blocks are more expensive to update than others for they are larger, bear more nonzero entries in the training set, suffer poorer data locality, or many other reasons. Blocks that take longer to update should have greater delays when they are updated because more other updates will have happened between the time that $\hat{x}$ is read and when the update is completed. For the same reason, when blocks are assigned agents, the updates by slower or busier agents will generally have greater delays.

Indeed this turns out to be the case. Experiments were performed on a cluster with 2 nodes, each with 16 threads running on an Intel Xeon CPU E5-2690 v2. The algorithm was applied to the logistic regression problem on the ``news20'' from LIBSVM, with 64 contiguous coordinate blocks of equal size. Over 2000 epochs, blocks 0, 1, and 15 have average delays of 351, 115, and 28,  respectively. %ASync-BCD completed this over 7x faster than the corresponding synchronous algorithm using the same computing resources, with a nearly equal decrease in objective function.

Even on a problem with each block having the same number of nonzeros, and when the computing environment is homogeneous, this dependence persists. We assigned 20 threads to each core, with each thread assigned to a block of 40 coordinates with equal numbers of nonzeros. The mean delay varied from 29 to 50 over the threads. This may be due to the cluster scheduler or issues of data locality, which were hard to examine. %More dramatic dependence was observed in a heterogeneous setting.

Clearly, there is strong dependence of the delays $\vec{j}(k)$ on the updated block $i_k$. Let us consider  an ideal situation where all blocks are equally difficult, all agents are equally fast, the job scheduling is fair, and data are centrally stored with the same distance to each agent. But, $\vec{j}(k)$ and $i_k$ still depend on the start times of the 1st through $k$th updates, so  $\vec{j}(k)$ and $i_k$  are still related. Therefore, it is  necessary to relax the independence assumption when applying the theory to asynchronous solvers.

%%%%%%%%%%%%%%%%%%%%%%%%%%%%%%%%%%%%%%%%
\subsection{Stochastic and deterministic block rules}\label{ssec:Setup}
%%%%%%%%%%%%%%%%%%%%%%%%%%%%%%%%%%%%%%%%
This paper considers two different \emph{block rules}: \emph{deterministic} and \emph{stochastic}.

For the stochastic block rule, for each update a block is chosen from $\{1,2,\ldots,N\}$ uniformly at random\footnote{The distribution doesn't have to be uniform. We need only assume that every block has a nonzero probability of being updated. It is easy to adjust our analysis to this case.}, for instance in \cite{LiuWrightReBittorfSridhar2015_asynchronous,LiuWright2015_asynchronous,peng2016arock}.
%$^,$\footnote{A stochastic rule is easier to analyze in general.}
For the deterministic rule, $i_k$ is an arbitrary sequence that is assumed to be \emph{essentially cyclic}. That is, there is an $N'\in\mathbb{N}$, $N'\ge N$, such that each block $i\in\{1,2,\ldots,N\}$ is updated at least once in a window of $N'$, that is,
\begin{center}For each $t\in \mathbb{Z}^+$, $\exists$ integer $K(i,t)\in\{tN',tN'+1,\ldots,(1+t)N'-1\}$ such that $i_{K(i,t)}=i$.
\end{center}
%That is, every block is updated at least once in every $N'$-sized window of iterations.
The essentially cyclic rule allows a cycle to go longer than $N$ due to update delays. It encompasses different kinds of \emph{cyclic} rules such as fixed ordering, random permutation, and greedy selection. %The ordering of each cycle is arbitrary. % in this paper is arbitrary.  %which may model many different implementation settings.

The stochastic block rule is easier to analyze because taking expectation will yield a good approximation to the full gradient. It ensures the every block is updated at the specified frequency. However, it can be expensive or even infeasible to implement for the following reasons.

In the shared memory setup, stochastic block rules require random data access, which is not only significantly slower than sequential data access but also cause frequent \emph{cache misses} (waiting for data being fetching from slower cache or the main memory).  The cyclic rules clearly avoid these issues.

In the server-worker setup where workers update randomly assigned blocks at each step, each worker must either store all the data or read the required data from the server at every step (in addition to reading $x$). This overhead is big. Permanently assigning blocks to agents is a more sound choice.

On the other hand, the analysis of cyclic rules generally has to consider the worst ordering and gives worse performance guarantees. %To guarantee convergence for any cyclic ordering, one must use a smaller step size, which leads to a slower rate.
%Nonetheless, these results hold tightly on only some problems and specific orderings of blocks.
In practice, cyclic rules often lead to good performance \cite{friedman2007pathwise,friedman2010regularization,ChowWuYin2016_cyclic}.

%The \emph{delayed iterate} $\widehat{x}^k$ is  what the updating agent has read from global memory to its local cache to compute its next update. We define $\vec{j}(k)=(j(k,1),j(k,2),\ldots,j(k,N))\in \mathbb{Z}^N$ as the \emph{delay vector} at the $k$-th iteration. The components of this vector represent how many iterations out-of-date the vector $\widehat{x}^k$ is\footnote{We follow the notation and definitions from \cite{hannah2016unbounded}.}. More precisely,
%
% \begin{align}
%     \hat{x}^k=(x^{k-j(k,1)}_1,x^{k-j(k,2)}_2,\ldots,x^{k-j(k,N)}_N).
% \end{align}
% %
% The \emph{$k$'th delay} (or \emph{current delay}) $j(k)\in \mathbb{N}$ is defined as $j(k)=\max_{1\leq l\leq N}\{j(k,l)\}$.

\subsection{Bounded and unbounded delays}
We consider different delay assumptions as well. \emph{Bounded delay} is when $j(k)\leq\tau$ for some fixed $\tau\in \mathbb{Z}^+$ and all iterations $k$; while the \emph{unbounded delay} allows $\sup_{k}\{j(k)\}=+\infty$. Bounded and unbounded delays can be further divided into deterministic and stochastic. Deterministic delays refer to the sequence of delay vectors $\vec{j}(0),\vec{j}(1),\vec{j}(2),\ldots$ that is arbitrary or follows an unknown distribution so is treated as arbitrary. Our stochastic delay results apply to distributions that decay faster than $O(k^{-3})$.

Deterministic unbounded delays apply to the case when async-BCD runs on unfamiliar hardware platforms. For convergence, we require a finite $\liminf_k \{j(k)\}$ and the current step size $\eta^k$ to be chosen adaptive to the current delay $j(k)$, which must be measured.

Bounded delays and stochastic unbounded delays apply when the user can provide the bound and the distribution, respectively. The user can obtain them from previous experience or by running a pilot test. In return, a fixed step size works, and measuring the current delay is not needed.

%%%%%%%%%%%%%%%%%%%%%%%%%%%%%%%%%%%%%%%%
\subsection{Contributions}
%%%%%%%%%%%%%%%%%%%%%%%%%%%%%%%%%%%%%%%%
The contributions are mainly convergence results for three kinds of delays: bounded, stochastic unbounded, deterministic unbounded, while allowing delays to depend on blocks. The results are provided for nonconvex, convex, and strongly convex functions with Lipschitz gradients. Sublinear rates and linear rates are provided, and, in terms of order of magnitude, they match their synchronous results.
%For bounded delays, we obtain convergence results (under constant step size) under both the deterministic  and stochastic block-index rules without the independence assumption, even when the objective function is nonconvex. A sublinear convergence rate can be derived if the function is convex, and a linear convergence rate is proved if the function is also restricted strongly convex for the stochastic index rule. We also prove a selection of other convergence results on unbounded delays. Other convergence results could easily be obtained using our approach given more space, but
Due to space limitation, we restrict ourselves to  Lipschitz differentiable functions and leave out nonsmooth proximable functions.

Like many analyses of algorithms, our proofs are built on the construction of Lyapunov functions. We provide a simple ODE-based (i.e., continuous time) construction for bounded delays. Once going discrete and considering the three different kinds of delays, the Lyapunov functions inevitably involve complicated coefficients. But, the ODE-based construction illustrates our construction principle and provides insights on how delays affect convergence.
%\commwy{Use rebuttal.} We are able to prove these results using technical innovations introduced in \cite{hannah2016unbounded} and  \cite{peng2016arock}. The use of Lyapunov functions allows us to prove much stronger results, by making ``sufficient descent'' results possible (e.g. Lemma \ref{lemma:bounded-descent}). Instead of analyzing, say, $f(x^k)$ or $\|x^k\|^2_2$, we add an ``asynchronicity error'' term, and analyze that instead (see (\ref{Lyapunov1})). Lyapunov functions appear to be the most natural error to consider when analyzing asynchronous algorithms. Also,

Roughly speaking, if the delays are bounded, then convergence should be analyzed over a sliding window of consecutive iterations because all the delays and progress in each window are intimately related. If no uniform bound is known, then the window must extend to the very first iteration.%, and the window size increases as the iteration does.
This analysis does bring great news to the practitioner. Basically speaking, even when there is no known load balancing (thus the delays may be sensible to the blocks) or bound of the delays, they can still ensure convergence by our provided step sizes. This applies to both random and deterministic choices of blocks.

We do not treat asynchronicity as noise as some recent papers do\footnote{See, for example, (5.1) and (A.10) in \cite{RechtReWrightNiu2011_hogwild}, (3.5) in \cite{mania2015perturbed}, and (14) and Lemma 4 in \cite{de2015taming}.}. In our setting, modelling asynchronicity in this way destroys valuable information, and leads to inequalities that are too blunt to obtain stronger results. Compared to noise, delays are much less harmful. This is why sublinear and linear rates can be established for weak and strong convex problems respectively, even when delays depend on the blocks and are potentially unbounded.
 %For instance, it can be shown that $\sum^\infty_{k=1}\|\nabla f(x^k)-\nabla f(\hat{x}^k)\|^2_2<\infty$ in some circumstances (e.g. using Lemma \ref{lemma:bounded-descent}), whereas this may not be true for a generic noise model.

We understand from the practitioner's point of view: there is a need for information on how to select parameters. However proving convergence in a new setting where there are no comparable convergence results at all takes a lot of the space. This prevents us from obtaining best-possible constants in our convergence rates, which can be extremely time-consuming. The main message is that sublinear and linear convergence can be obtained under more realistic assumptions. %in the appropriate settings. Better rates are possible by optimizing the coefficients of the Lyapunov function, but
Our results do provide step sizes when the delay bound is given, the delay distribution is known, or the current delay is measured. The convergence results are generally tight in terms of the order of the involved quantities, but the constants are perhaps not tight .

%%%%%%%%%%%%%%%%%%%%%%%%%%%%%%%%%%%%%%%
\subsection{Related work}
%%%%%%%%%%%%%%%%%%%%%%%%%%%%%%%%%%%%%%%
Our work extends the theory on asynchronous BCD algorithms such as \cite{RechtReWrightNiu2011_hogwild,LiuWrightReBittorfSridhar2015_asynchronous,LiuWright2015_asynchronous}. However, their analysis relies the independence assumption and assume bounded delays. The bounded delay assumption was weakened by recent papers \cite{hannah2016unbounded,pengxuyanyin17}, but independence and random blocks were still needed. %, and \cite{hannah2016unbounded} focuses the fixed-point problem without an objective function and allows only the random block rule.
%\rev{Many existing analyses for  treats delays as noise and manage to bound it by the progress at each step. (Examples include Eqs. (5.1) and (A.10) in \cite{RechtReWrightNiu2011_hogwild}; Eq.  (3.5) in \cite{mania2015perturbed}; (14)  and Lemma 4 in \cite{de2015taming}.)}

Recently \cite{LeblondPedregosaLacoste-Julien2016_asaga} proposes (in the SGD setting) %that the independence assumption appeared in almost of the relevant literature, and that enforcing it would require synchronization or locks. They
a novel ``read after'' sequence relabeling technique to create the independence. % between $i_k$ and the delay.
However, enforcing independence in this way creates other artificial implementation requirements that may waste computational resources: For instance, agents need to read \emph{all} shared data before computing their update, %forcing computation to stop.  The entire shared vector must be read
even if not all of it is required to compute updates, which can be extremely expensive for sparse data. It is also necessary to recompute certain parameters to prevent a biased update estimator, instead of caching and cheaply updating. Our analysis does not require these kinds of implementation fixes because it does not rely on any kind of unbiased update estimator. %Reading and computing can happen simultaneously.
Also, our analysis also works for unbounded delays and deterministic block choices.

Related recent works also include \cite{cannelli2016asynchronous,cannelli2017asynchronous}, which solve our problem with additional convex block-separable terms in the objective. In the first paper \cite{cannelli2016asynchronous}, independence between blocks and delays is avoided. However, they require a step size that diminishes at $1/k$ and that the sequence of iterate is bounded (which in general may not be true). The second paper \cite{cannelli2017asynchronous} relaxes independence by using a different set of assumptions. In particular, their assumption D3 assumes that, regardless of the previous updates, there is a universally positive chance for every block to be updated in the next step. This Markov-type assumption relaxes the independence assumption but does not avoid it. In particular, it is not satisfied by Example 1 below. %the small example below \eqref{eq:condEcond}. %Assumption D4 assigns each block to an agent permanently, thus allowing \cite{cannelli2017asynchronous}  to use proximal mapping for the block-separable nonsmooth terms in their objective. Like \cite{LiuWright2015_asynchronous,davis2016asynchronous}, proximal mappings must use their up-to-date blocks of variables, which is ensured by Assumption D4 in \cite{cannelli2017asynchronous}. So far, only operator-based approaches \cite{peng2016arock,hannah2016unbounded} allow proximal mappings to use delayed blocks. % \rev{the authors still admit the independent assumption.}
%In the nonconvex case, the step size in \cite{cannelli2016asynchronous,cannelli2017asynchronous} must be decreasing while ours is fixed.

In the convex case with a bounded delay $\tau$, the step size in paper \cite{LiuWrightReBittorfSridhar2015_asynchronous} is
{$O(\frac{1}{\tau^2/N})$}. In their proofs, the Lyapunov function is based on $\|x^{k}-x^*\|_2^2$.  Our analysis uses a Lyapunov function consisting of both the function value and the sequence history, where the latter vanishes when delays vanish. If the $\tau$ is much larger than the blocks of the problem,
our result {$O(\frac{1}{\tau})$} is better even under our much weaker conditions.  The step size bound  in \cite{peng2016arock,hannah2016unbounded,davis2016asynchronous} is {$O(\frac{1}{{\tau}/{\sqrt{N}}})$}, which is better than ours, but they need the independence assumption and the stochastic block rule.

%The paper \cite{davis2016asynchronous} uses the constant step size $\frac{1}{O(\frac{\tau}{\sqrt{N}})}$ like in  \cite{peng2016arock,hannah2016unbounded}.

Recently, \cite{xu2017asynchronous} introduces an asynchronous method primal-dual for a problem similar to ours but having additional affine linear constraints. The analysis assumes bounded delays, random blocks, and independence.

%%%%%%%%%%%%%%%%%%%%%%%%%%%%%%%%%%%%%%%
\subsection{Notation}
%%%%%%%%%%%%%%%%%%%%%%%%%%%%%%%%%%%%%%%
Let $x^*\in \textrm{arg}\min f$. For the update in \eqref{eq:ARock-algorithm}, we use the following notation:
\small
\begin{align}\label{eq:Dddef}
\Delta^k:=x^{k+1}-x^k\overset{\eqref{eq:ARock-algorithm}}{=}-\frac{\gamma_k}{L}\nabla_{i_k},\qquad d^k:=x^k-\hat{x}^k.
\end{align}
\normalsize
We also use the convention $\Delta^k:=0$ if $k<0$.
Let $\chi^k$ be the sigma algebra generated by $\{x^0,x^1,\ldots,x^k\}$. %i.e., $\chi^k=\sigma(x^0,x^1,\ldots,x^k)$.
Let $\EE_{\vec{j}(k)}$ denote the expectation over the value of $\vec{j}(k)$ (when it is a random variable). $\EE$ denotes the total expectation.

%%%%%%%%%%%%%%%%%%%%%%%%%%%%%%%%%%%%%%
\section{Bounded delays}
%%%%%%%%%%%%%%%%%%%%%%%%%%%%%%%%%%%%%%%%
In this part, we present convergence results for the bounded delays. If the gradient of the function is $L$-Lipschitz (even if the function is nonconvex), we present the convergence for both the deterministic and stochastic  block rule. If the function is convex, we can obtain a sublinear convergence rate.
Further, if the function is restricted strongly convex, a linear convergence rate is obtained.

\subsection{Continuous-time analysis}
%We present a continuous version of our Lyapunov function.
Let $t$ be time in this subsection. Consider the ODE
\small
\begin{align}\label{eq:delayedgradflow}
    \dot{x}(t)=-\eta\nabla f(\hat{x}(t)),
\end{align}
\normalsize
where $\eta >0$.
If we set $\hat{x}(t)\equiv x(t)$, this system describes a gradient flow, which monotonically decreases $f(x(t))$, and  its discretization is the gradient descent iteration. Indeed, we have
$$
\frac{d}{dt} f(x(t)) = \langle\nabla f(x(t)),\dot{x}(t)\rangle \overset{\eqref{eq:delayedgradflow}}{=} - \frac{1}{\eta}\|\dot x(t)\|_2^2.
$$

Instead, we allow delays (i.e., $\hat{x}(t)\neq x(t)$) and impose the bound $c>0$ on the delays:
\small
\begin{align}
\label{eq:debnd}
\|\hat{x}(t)-x(t)\|_2\leq \int_{t-c}^{t}\|\dot{x}(s)\|_2 ds.
\end{align}
\normalsize
The delays introduce inexactness to the gradient flow $f(x(t))$. We lose monotonicity. Indeed,
\small
\begin{align}
    &\frac{d}{dt} f(x(t))=\langle\nabla f(x(t)),\dot{x}(t)\rangle=\langle\nabla f(\hat{x}(t)),\dot{x}(t)\rangle+\langle\nabla f(x(t))-\nabla f(\hat{x}(t)),\dot{x}(t)\rangle\label{eq:contdotp}\\
    &\quad\overset{a)}{\leq}-\frac{1}{\eta}\|\dot{x}(t)\|_2^2+L\|x(t)-\hat{x}(t)\|_2\cdot\|\dot{x}(t)\|_2\overset{b)}{\leq} -\frac{1}{2\eta}\|\dot{x}(t)\|_2^2+\frac{\eta c L^2}{2}\int_{t-c}^{t}\|\dot{x}(s)\|_2^2ds,\nonumber
\end{align}
\normalsize
where a) is from \eqref{eq:delayedgradflow} and Lipschitzness of $\nabla f$ and b) is from the Cauchy-Schwarz inequality $L\|x(t)-\hat{x}(t)\|_2\cdot\|\dot{x}(t)\|_2\leq \frac{\|\dot{x}(t)\|_2^2}{2\eta}+\frac{\eta L^2\|x(t)-\hat{x}(t)\|_2^2}{2}$ and, by \eqref{eq:debnd}, $\|x(t)-\hat{x}(t)\|_2^2\leq c\int^t_{t-c}\|\dot{x}(s)\|_2^2ds$.
The inequalities are generally unavoidable.
%without more properties of $f$.
The integral term is due to the use of delayed gradient.

Therefore, we design an energy with both $f$ and a weighted total kinetic term, where $\gamma >0$ will be decided below:
\small
\begin{align}
    \xi(t)=f(x(t))+\gamma\int_{t-c}^{t}\big(s-(t-c)\big)\|\dot{x}(s)\|^2_2ds.
\end{align}
\normalsize
%We argue that the delayed gradient flow \eqref{eq:delayedgradflow} reduces this energy monotonically. Indeed, with the integration by part,
$\xi(t)$ has the time derivative
\small
\begin{align*}
    \dot{\xi}(t)&=\frac{d}{dt} f(x(t))+\gamma c\|\dot{x}(t)\|_2^2-\gamma\int_{t-c}^{t}\|\dot{x}(s)\|_2^2ds.
\end{align*}
\normalsize
By substituting the bound on $\frac{d}{dt}f(x(t))$ in \eqref{eq:contdotp}, we get
\small
\begin{align}\label{eq:dotxibnd}
    \dot{\xi}(t)&\le -(\frac{1}{2\eta}-\gamma)\|\dot{x}(t)\|_2^2-(\gamma-\frac{\eta c L^2}{2})\int_{t-c}^{t}\|\dot{x}(s)\|_2^2ds.
\end{align}
\normalsize
As long as $\eta<\tfrac{1}{Lc}$, there exists $\gamma>0$ such that $(\frac{1}{2\eta}-\gamma)> 0$ and $(\gamma-\frac{\eta c L^2}{2})> 0$, so $\xi(t)$ is monotonically nonincreasing. %which together ensure the monotonicity of the energy $\xi(t)$.
Assume $\min f$ is finite. Since $\xi(t)$ is lower bounded by $\min f$, $\xi(t)$ must converge, subsequently yielding the convergence of $\dot{\xi}\to 0$, $\dot{x}(t)\to 0$ by \eqref{eq:dotxibnd}, $\nabla f(\hat{x}(t))\to 0$ by \eqref{eq:delayedgradflow},  and $\hat{x}(t)-x(t)\to 0$ by \eqref{eq:debnd}. The last two results further give $\nabla f(x(t))\to 0$.

\subsection{Discrete analysis}
The analysis for our discrete iteration \eqref{eq:ARock-algorithm} is based on the following Lyapunov function:
\small
\begin{align}
\xi_k:=f(x^k)+\frac{L}{2\varepsilon}\sum_{i=k-\tau}^{k-1}(i-(k-\tau)+1)\|\Delta^i\|_2^2.\label{Lyapunov1}
\end{align}
\normalsize
for $\varepsilon>0$ to be determined later based on the step size and $\tau$, the bound on the delays.

In the lemma below, we present a fundamental inequality, which states, regardless of which block $i_k$ is updated and which $\hat{x}^k$ is used to compute the update in \eqref{eq:ARock-algorithm}, there is a sufficient descent in our Lyapunov function for a proper step size  $1/O(\tau)$. %, convergence of the updates $\Delta^k$ is ensured.
%%%%%%%%%%%%%%%%%%%%%%%%%%%%%%%%%%%%%%%%
\begin{lemma}[sufficient descent for bounded delays] \label{lemma:bounded-descent}
%%%%%%%%%%%%%%%%%%%%%%%%%%%%%%%%%%%%%%%%
\textbf{Conditions:} Let $f$ be a function (possibly nonconvex) with $L$-Lipschitz gradient and finite $\min f$. Let $(x^k)_{k\geq 0}$ be generated by the async-BCD algorithm \eqref{eq:ARock-algorithm}, and the delays be bounded by $\tau$. Choose the step size $$\gamma_k\equiv\gamma=\frac{2c}{2\tau+1}$$ for arbitrary fixed $0<c<1$.
\textbf{Result:} we can choose $\varepsilon>0$ to obtain
\small
\begin{align}\label{bounded-descent-result-1}
    \xi_k-\xi_{k+1}\geq\frac{1}{2}(\frac{1}{\gamma}-\frac{1}{2}-\tau)L\cdot\|\Delta^k\|_2^2,
\end{align}
\normalsize
Consequently,
\vspace{-10pt}
\small
\begin{align}
 \lim_{k}\|\Delta^k\|_2&=0\label{bounded-descent-result-2},\\
 \min_{1\leq i\leq k}\|\Delta^i\|_2&=o(1/\sqrt{k}).\label{bounded-descent-result-3}
\end{align}
\normalsize
\end{lemma}
The rate in \eqref{bounded-descent-result-3} concerns the smallest $\|\Delta^i\|_2$ among $i=1,\ldots,k$. We call it the \emph{running best rate}. Although $\|\Delta^k\|_2$ is not monotonic, \eqref{bounded-descent-result-2} and \eqref{bounded-descent-result-3} indicate that $\|\Delta^k\|_2$ decays \emph{overall faster} than $1/\sqrt{k}$. %Also, we keep the square in \eqref{bounded-descent-result-3} since it naturally appears in the analysis, and the literature often (but not always) uses it with the square.
%%%%%%%%%%%%%%%%%%%%%%%%%%%%%%%%%%%%%%%%

%
%%%%%%%%%%%%%%%%%%%%%%%%%%%%%%%%%%%%%%%%

%%%%%%%%%%%%%%%%%%%%%%%%%%%%%%%%%%%%%%%%
%\subsection{Nonconvex, with deterministic block rule}
%%%%%%%%%%%%%%%%%%%%%%%%%%%%%%%%%%%%%%%%
Based on the lemma, we obtain a very general result for nonconvex problems: % and the deterministic block rule as follows
%%%%%%%%%%%%%%%%%%%%%%%%%%%%%%%%%%%%%%%%
\begin{theorem}\label{thm:bounded-deterministic}
%%%%%%%%%%%%%%%%%%%%%%%%%%%%%%%%%%%%%%%%
Assume the conditions of Lemma \ref{lemma:bounded-descent}, for $f$ that may be nonconvex. Under the deterministic block rule, we have
\small
\begin{align}\label{bounded-deterministic-result-1}
    \lim_{k}\|\nabla f(x^k)\|_2=0,\quad \min_{1\le i\le k}\|\nabla f(x^k)\|_2 = o(1/\sqrt{k}).
\end{align}
\normalsize
\end{theorem}
This rate has the same order of magnitude as standard gradient descent.
%%%%%%%%%%%%%%%%%%%%%%%%%%%%%%%%%%%%%%%%

%%%%%%%%%%%%%%%%%%%%%%%%%%%%%%%%%%%%%%%%
\subsection{Stochastic block rule}
%%%%%%%%%%%%%%%%%%%%%%%%%%%%%%%%%%%%%%%%
Under the stochastic block rule, an agent picks a block from $\{1,2,...,N\}$ uniformly randomly at the beginning of each update. For the $k$th completed update, the index of the chosen block is $i_k$.
Our result in this subsection relies on the following assumption on the random variable $i_k$:
\begin{align}\label{eq:condEcond}
 \EE_{i_k}(\|\nabla_{i_k} f(x^{k-\tau})\|_2\mid\chi^{k-\tau})=\frac{1}{N}\sum_{i=1}^N\|\nabla_{i} f(x^{k-\tau})\|_2,
\end{align}
where $\chi^k=\sigma(x^0,x^1,\ldots,x^k,\vec{j}(0),\vec{j}(1),\ldots,\vec{j}(k))$, $k=0,1,\ldots$, is the filtration that represents the information that is accumulated as our algorithm runs. It is important to note that \eqref{eq:condEcond} uses $x^{k-\tau}$ instead of $\hat{x}^k$. The latter is \emph{not} independent of $i_k$.

Condition \eqref{eq:condEcond} is a property about $i_k$ rather than $\nabla f$. It states that each of the $N$ possible values of $i_k$ occurs at probability $1/N$ given the information $\tau$ and more iterations older. \emph{We can relax \eqref{eq:condEcond} to non-uniform distributions}; indeed, Theorem \ref{thm:bounded-stochastic-general} below only needs that every block has a nonzero probability of being updated given $\chi^{k-\tau}$, that is,
\begin{align}\label{eq:condEcond11}
 \EE(\|\nabla_{i_k} f(x^{k-\tau})\|_2\mid\chi^{k-\tau})\ge\frac{\varepsilon}{N}\sum_{i=1}^N\|\nabla_{i} f(x^{k-\tau})\|_2,
\end{align}
for some universal $\varepsilon>0$.
The uniform distribution in Assumption \eqref{eq:condEcond} is made for convenience and simplicity.

The above assumption is justified since this section assumes a finite delay bound $\tau$ and thus the history older than $\tau$, though might still affect $i_k$, can no longer nullify the chance of $i_k$ taking each of $\{1,\ldots,N\}$. On the other hand, making a similar assumption on $\EE(\|\nabla_{i_k} f(x^{k-\tau})\|_2\mid\chi^{t})$, for any $t=k-1,k-2\ldots,k-\tau+1$, would be \emph{unjustified}, as shown in the following example.

\noindent\textbf{Example 1.} Consider three different blocks and two identical agents. Assume blocks 1,2,3 take exactly 2,3,4 seconds to update by either agent. The maximal delay is $\tau=4/2=2$. %At each step, each agent picks one of three blocks at random.
Assume both agents start their first jobs at nearly the same time. If the first completed update is $i_1=2$, by one of two agents, then $i_2$ must equal either $2$ or $3$; $i_2=1$ is impossible. This can be verified by enumerating all the possible combinations of the block choices made by the agents in their first two steps. In general, $i_k=1$ is impossible when, before the $(k-1)$th completed update, the two agents start their new steps at nearly the same time and $i_{k-1}=2$.

%%%%%%%%%%%%%%%%%%%%%%%%%%%%%%%%%%%%%%%%
%\subsubsection{General nonconvex result}
%%%%%%%%%%%%%%%%%%%%%%%%%%%%%%%%%%%%%%%%

Next, we present a general result for a possibly nonconvex objective $f$.
%%%%%%%%%%%%%%%%%%%%%%%%%%%%%%%%%%%%%%%%
\begin{theorem}\label{thm:bounded-stochastic-general}
%%%%%%%%%%%%%%%%%%%%%%%%%%%%%%%%%%%%%%%%
Assume the conditions of Lemma \ref{lemma:bounded-descent}.%, for $f$ that may be nonconvex.
Under the stochastic block rule and assumption \eqref{eq:condEcond}, we have:
\small
\begin{align}\label{bounded-stochastic-general-result-1}
    \lim_{k}\EE\|\nabla f(x^k)\|_2=0,\quad \min_{1\le i\le k}\EE\|\nabla f(x^k)\|_2^2 = o(1/k).
\end{align}
\normalsize
\end{theorem}
%%%%%%%%%%%%%%%%%%%%%%%%%%%%%%%%%%%%%%%%

%%%%%%%%%%%%%%%%%%%%%%%%%%%%%%%%%%%%%%%%

%%%%%%%%%%%%%%%%%%%%%%%%%%%%%%%%%%%%%%%%

%%%%%%%%%%%%%%%%%%%%%%%%%%%%%%%%%%%%%%%%
%\begin{remark}
%%%%%%%%%%%%%%%%%%%%%%%%%%%%%%%%%%%%%%%%
%Under the independence assumption, the stepsize bound is $\frac{2}{\frac{2\tau}{\sqrt{N}}+1}$ \cite{peng2016arock}. Without the independence assumption, the stepsize bound is smaller.
%
%\end{remark}
%%%%%%%%%%%%%%%%%%%%%%%%%%%%%%%%%%%%%%%%

%%%%%%%%%%%%%%%%%%%%%%%%%%%%%%%%%%%%%%%%
\subsubsection{Sublinear rate under convexity}
%%%%%%%%%%%%%%%%%%%%%%%%%%%%%%%%%%%%%%%%
When the function $f$ is convex, we can obtain convergence rates, for which we need a slightly modified Lyapunov function
\small
\begin{align}\label{Lyapunov2}
   F_k:=f(x^k)+\delta\cdot\sum_{i=k-\tau}^{k-1}(i-(k-\tau)+1)\|\Delta^i\|_2^2,
\end{align}
\normalsize
where $\delta:=[1+\frac{\varepsilon}{2\tau}(\frac{1}{\gamma}-\frac{1}{2}-\tau)]\frac{L}{2\varepsilon}$\footnote{Here, we assume $\tau\geq 1$.}.
We also define
$\pi_k:=\EE(F_k-\min f),\quad S(k,\tau):=\sum_{i=k-\tau}^{k-1}\delta\|\Delta^i\|_2^2.$
%to present the following lemma.

%%%%%%%%%%%%%%%%%%%%%%%%%%%%%%%%%%%%%%%%
\begin{lemma} \label{lem:bounded-stochastic-convex}
%%%%%%%%%%%%%%%%%%%%%%%%%%%%%%%%%%%%%%%%
Assume the conditions of Lemma \ref{lemma:bounded-descent}. Furthermore, let $f$ be convex and use the stochastic block rule. Let $\overline{x^k}$ denote the projection of $x^k$ to $\argmin{f}$, assumed to exist, and let
\small
\begin{align}\label{eq:alpha}
\beta :=\max\{\tfrac{8NL^2}{\gamma^2},(12N+2)L^2\tau+\delta\tau\}, \quad \alpha: =\beta/[\tfrac{L}{4\tau}(\tfrac{1}{\gamma}-\tfrac{1}{2}-\tau)].
\end{align}
\normalsize
Then we have:
\small
\begin{align}\label{bounded-stochastic-convex-result}
   (\pi_{k})^2\leq\alpha(\pi_{k}-\pi_{k+1})\cdot(\tau\EE S(k,\tau)+\EE\|x^k-\overline{x^k}\|_2^2).
\end{align}
\normalsize
\end{lemma}
%%%%%%%%%%%%%%%%%%%%%%%%%%%%%%%%%%%%%%%%
When $\tau=1$ (nearly no delay), we can obtain $\beta=O(NL^2/\gamma^2)$ and $\alpha = O(\beta\gamma /L) = O(NL/\gamma)$, which matches the result of standard BCD. Unfortunately, delays cause the complication of $\alpha$.
%%%%%%%%%%%%%%%%%%%%%%%%%%%%%%%%%%%%%%%%

%%%%%%%%%%%%%%%%%%%%%%%%%%%%%%%%%%%%%%%%

We now present the sublinear convergence rate.

%%%%%%%%%%%%%%%%%%%%%%%%%%%%%%%%%%%%%%%%
\begin{theorem}\label{thm:bounded-stochastic-convex}
%%%%%%%%%%%%%%%%%%%%%%%%%%%%%%%%%%%%%%%%
Assume the conditions of Lemma \ref{lemma:bounded-descent}. Furthermore, let $f$ be convex and coercive\footnote{A function $f$ is coercive if  $\|x\|\to\infty$ means $f(x)\to\infty$.}, and use the stochastic block rule.
Then we have:
\small
\begin{align}
    \EE(f(x^k)-\min f)= O(1/k).
\end{align}
\normalsize
\end{theorem}
%%%%%%%%%%%%%%%%%%%%%%%%%%%%%%%%%%%%%%%%

%%%%%%%%%%%%%%%%%%%%%%%%%%%%%%%%%%%%%%%%

%%%%%%%%%%%%%%%%%%%%%%%%%%%%%%%%%%%%%%%%

%%%%%%%%%%%%%%%%%%%%%%%%%%%%%%%%%%%%%%%%
\subsubsection{Linear rate under convexity}
%%%%%%%%%%%%%%%%%%%%%%%%%%%%%%%%%%%%%%%%
We next consider when $f$ is $\nu$-\textit{restricted strongly convex}\footnote{A condition weaker than $\nu$-strong convexity and useful for problems involving an underdetermined linear mapping $Ax$; see \cite{lai2013augmented,LiuWright2015_asynchronous}.} in addition to having $L$-Lipschitz gradient. That is, for $x\in \textrm{dom}(f)$,
\small
\begin{align}\label{RSC}
    \langle\nabla f(x),x-\textrm{Proj}_{\argmin{f}}(x)\rangle\geq \nu\cdot\textrm{dist}^2(x,\argmin{f}).
\end{align}
\normalsize
%

%%%%%%%%%%%%%%%%%%%%%%%%%%%%%%%%%%%%%%%%
\begin{theorem} \label{thm:bounded-stochastic-restrictedstronglyconvex}
%%%%%%%%%%%%%%%%%%%%%%%%%%%%%%%%%%%%%%%%
Assume the conditions of Lemma \ref{lemma:bounded-descent}. Furthermore, let $f$ be $\nu$-strongly convex, and use the stochastic block rule. Then we have:
\small
\begin{align}
    \EE(f(x^k)-\min f)= O(c^{k}),
\end{align}
\normalsize
where $c:={\tfrac{\alpha }{\min\{\nu,1\}}}\big/{(1+\tfrac{\alpha }{\min\{\nu,1\}})}<1$ for $\alpha$ given in \eqref{eq:alpha}.

\end{theorem}
%%%%%%%%%%%%%%%%%%%%%%%%%%%%%%%%%%%%%%%%

%%%%%%%%%%%%%%%%%%%%%%%%%%%%%%%%%%%%%%%%

%%%%%%%%%%%%%%%%%%%%%%%%%%%%%%%%%%%%%%
\section{Stochastic unbounded delay}
%%%%%%%%%%%%%%%%%%%%%%%%%%%%%%%%%%%%%%%%
%The major difference from previous proofs for bounded case lies in dealing with the term $\EE\|d^{k}\|_2^2$.
In this part, the delay vector $\vec{j}(k)$ is an unbounded random variable, which allow extremely large delays in our algorithm. Under some mild restrictions on the distribution of $\vec{j}(k)$, we can still establish convergence. In light of our continuous-time analysis, we must develop a new bound for the last inner product in \eqref{eq:contdotp}, which requires the tail distribution of $j(k)$ to decay sufficiently fast.

%\rev{Verify whether the distribution can be changing. Probably is. Then emphasise this: We don't need a fixed distribution. Very general etc.}

Specifically, we define fixed parameters $p_j$ such that $p_j\geq\mathbb{P}(j(k)=j), \forall k$, $s_l=\sum_{j=l}^{+\infty} jp_j$, and $c_i:=\sum_{l=i}^{+\infty}s_l$. Clearly, $c_0$ is larger than $c_1,c_2,\ldots$, and we need $c_0$ to be finite.
Distributions with $p_j=\mathcal{O}(j^{-t})$, $t>4$, and exponential-decay distributions satisfy this requirement.

Define the Lyapunov function $G_k$ as
\small
\begin{align}\label{Lyapunov3}
   G_k:=f(x^k)+\bar{\delta}
    \cdot\sum_{i=0}^{k-1}c_{k-1-i}\|\Delta^{i}\|_2^2,
\end{align}
\normalsize
where $\bar{\delta}:=\tfrac{L}{2\varepsilon}+(\tfrac{1}{\gamma}-\tfrac{1}{2})\tfrac{L}{c_0}-\frac{L}{\sqrt{c_0}}$. To simplify the presentation, we define $R(k):=\sum_{i=0}^{k}c_{k-i}\EE\|\Delta^{i}\|_2^2.$
%
%%%%%%%%%%%%%%%%%%%%%%%%%%%%%%%%%%%%%%%%
\begin{lemma}[Sufficient descent for stochastic unbounded delays] \label{lemma:stochastic-unbounded-descent}
%%%%%%%%%%%%%%%%%%%%%%%%%%%%%%%%%%%%%%%%
\textbf{Conditions:} Let $f$ be a function (which may be nonconvex) with $L$-Lipschitz gradient and finite $\min f$. Let delays be stochastic unbounded. %Let $(x^k)_{k\geq 0}$ be generated by async-BCD.
Use step size $\gamma_k\equiv\gamma=\frac{2c}{2\sqrt{c_0}+1}$ for arbitrary fixed $0<c<1$. \textbf{Results:} we can set $\varepsilon>0$ to ensures sufficient descent:
\small
\begin{align}\label{stochastic-unbounded-descent-result-1}
\EE[G_k-G_{k+1}]\geq\tfrac{L}{c_0}(\tfrac{1}{\gamma}-\tfrac{1}{2}-\sqrt{c_0})R(k).
\end{align}
\normalsize
And we have
\small
\begin{align}\label{stochastic-unbounded-descent-result-2}
  \lim_{k}\EE\|\Delta^k\|_2=0~~\text{ and }~~\lim_{k}\EE\|d^k\|_2=0.
\end{align}
\normalsize
\end{lemma}
%%%%%%%%%%%%%%%%%%%%%%%%%%%%%%%%%%%%%%%%
For technical reasons, it appears to be difficult to obtain the rates of $\EE\|\Delta^k\|_2$ and $\EE\|d^k\|_2$.
%%%%%%%%%%%%%%%%%%%%%%%%%%%%%%%%%%%%%%%%

%%%%%%%%%%%%%%%%%%%%%%%%%%%%%%%%%%%%%%%%
\subsection{Deterministic block rule}
%%%%%%%%%%%%%%%%%%%%%%%%%%%%%%%%%%%%%%%%

%%%%%%%%%%%%%%%%%%%%%%%%%%%%%%%%%%%%%%%%
\begin{theorem}\label{thm:stochastic-unbounded-deterministic}
%%%%%%%%%%%%%%%%%%%%%%%%%%%%%%%%%%%%%%%%
Let the conditions of Lemma \ref{lemma:stochastic-unbounded-descent} hold for $f$. Under the deterministic block rule, we have:
\small
\begin{align}\label{stochastic-unbounded-deterministic-result-1}
    \lim_{k}\EE\|\nabla f(x^k)\|_2=0.
\end{align}
\normalsize
\end{theorem}
%%%%%%%%%%%%%%%%%%%%%%%%%%%%%%%%%%%%%%%%

%%%%%%%%%%%%%%%%%%%%%%%%%%%%%%%%%%%%%%%%

%%%%%%%%%%%%%%%%%%%%%%%%%%%%%%%%%%%%%%%%

%%%%%%%%%%%%%%%%%%%%%%%%%%%%%%%%%%%%%%%%
\subsection{Stochastic block rule}
%%%%%%%%%%%%%%%%%%%%%%%%%%%%%%%%%%%%%%%%

Recall that under the stochastic block rule, the block to update is selected uniformly at random from $\{1,2,\ldots,N\}$.
The previous assumption \eqref{eq:condEcond}, which is made for bounded delays, need to be updated into the following assumption for unbounded delays:
\begin{align}\label{unassump}
 \EE_{i_k}(\|\nabla_{i_k}f(x^{k-j(k)})\|_2^2)=\frac{1}{N}\sum_{i=1}^N\|\nabla_{i} f(x^{k-j(k)})\|_2^2,
\end{align}
where $j(k)$ is a random variable on both sides.
As argued below \eqref{eq:condEcond}, the uniform distribution can be relaxed, but we use it for simplicity.

%%%%%%%%%%%%%%%%%%%%%%%%%%%%%%%%%%%%%%%%
\begin{theorem} \label{thm:stochastic-unbounded-stochastic}
%%%%%%%%%%%%%%%%%%%%%%%%%%%%%%%%%%%%%%%%
Let the conditions of Lemma \ref{lemma:stochastic-unbounded-descent} hold. Under the stochastic block rule and assumption \eqref{unassump},  we have
\small
\begin{align}
    \lim_{k}\EE\|\nabla f(x^k)\|_2=0.
\end{align}
\normalsize
\end{theorem}
%%%%%%%%%%%%%%%%%%%%%%%%%%%%%%%%%%%%%%%%

%%%%%%%%%%%%%%%%%%%%%%%%%%%%%%%%%%%%%%%%

%%%%%%%%%%%%%%%%%%%%%%%%%%%%%%%%%%%%%%%%
\subsubsection{Convergence rate}
%%%%%%%%%%%%%%%%%%%%%%%%%%%%%%%%%%%%%%%%
When $f$ is convex, we can derive convergence rates for $$\phi_k:=\EE(G_k-\min f).$$
%%%%%%%%%%%%%%%%%%%%%%%%%%%%%%%%%%%%%%%%
\begin{lemma}\label{lem:stochastic-unbounded-stochastic-convex}
%%%%%%%%%%%%%%%%%%%%%%%%%%%%%%%%%%%%%%%%
Let the conditions of Lemma \ref{lemma:stochastic-unbounded-descent} hold, and let $f$ be convex. Let $\overline{x^k}$ denote the projection of $x^k$ to $\argmin f$. Let $\overline{\beta}=\max\{\frac{8NL^2}{\gamma^2c_0},(12N+2)L^2+\bar{\delta}\}$ and $\overline{\alpha}=\overline{\beta}/[\frac{L}{2}(\tfrac{1}{\gamma}-\tfrac{1}{2}-\sqrt{c_0})]$. Then we have
\small
\begin{align}
   (\phi_{k})^2\leq\bar{\alpha}(\phi_{k}-\phi_{k+1})\cdot(\bar{\delta}R(k)+\EE\|x^k-\overline{x^k}\|_2^2),
\end{align}
\normalsize
\end{lemma}
%%%%%%%%%%%%%%%%%%%%%%%%%%%%%%%%%%%%%%%%

%%%%%%%%%%%%%%%%%%%%%%%%%%%%%%%%%%%%%%%%
%%%%%%%%%%%%%%%%%%%%%%%%%%%%%%%%%%%%%%%%

A sublinear convergence rate can be obtained if $\sup_{k}\{\EE\|x^k-\overline{x^k}\|_2^2\}<+\infty$, which can be ensured by adding a projection to a large artificial box set that surely contains the solution. Here we only present a linear convergence result.

%%%%%%%%%%%%%%%%%%%%%%%%%%%%%%%%%%%%%%%%
\begin{theorem}\label{lem:stochastic-unbounded-stochastic-convex-linear}
%%%%%%%%%%%%%%%%%%%%%%%%%%%%%%%%%%%%%%%%
Let the conditions of Lemma \ref{lemma:stochastic-unbounded-descent} hold. In addition, let $f$ be $\nu$-restricted strongly convex and set step size $\gamma_k\equiv\gamma<\frac{2}{2\sqrt{c_0}+1}$, with $c=\frac{\bar{\alpha}\max\{1,\frac{1}{\nu}\}}{1+\bar{\alpha}\max\{1,\frac{1}{\nu}\}}<1$. Then,
\small
\begin{equation}
    \EE(f(x^k)-\min f)= O(c^{k}).
\end{equation}
\normalsize
\end{theorem}
%%%%%%%%%%%%%%%%%%%%%%%%%%%%%%%%%%%%%%%%

%%%%%%%%%%%%%%%%%%%%%%%%%%%%%%%%%%%%%%%%

%%%%%%%%%%%%%%%%%%%%%%%%%%%%%%%%%%%%%%%%

%%%%%%%%%%%%%%%%%%%%%%%%%%%%%%%%%%%%%%%%

%%%%%%%%%%%%%%%%%%%%%%%%%%%%%%%%%%%%%%%%

%%%%%%%%%%%%%%%%%%%%%%%%%%%%%%%%%%%%%%%%

%%%%%%%%%%%%%%%%%%%%%%%%%%%%%%%%%%%%%%
\section{Deterministic unbounded delays}
%%%%%%%%%%%%%%%%%%%%%%%%%%%%%%%%%%%%%%%%
In this part, we consider deterministic unbounded delays, which require delay-adaptive step sizes.
Set positive sequence $(\epsilon_i)_{i\geq 0}$ (which can be optimized later given the delays) such that $\kappa_i:=\sum_{j=i}^{+\infty}\epsilon_j$ obeys $\kappa_1<+\infty$. Set $D_j:=\frac{1}{2}+\tfrac{ \kappa_1}{2}+\sum_{i=1}^{j}\tfrac{1}{2\epsilon_i}$.
We use a new Lyapunov function:
%$H_k$ as
\small
\begin{eqnarray}\label{Lya4}
   H_k:=f(x^k)+\tfrac{L}{2}\sum_{i=1}^{+\infty}\kappa_{i}\|\Delta^{k-i}\|_2^2.
\end{eqnarray}
\normalsize
For any $T\geq\lim\inf{j(k)}$, let $Q_{T}$ be the subsequence of $\mathbb{N}$ where the current delay is less than $T$. The points $x^k,~k\in Q_T$, have convergence guarantees. The remaining points are still computed, but because they are affected by potentially unbounded delays, their quality is not guaranteed. % our algorithm should stop at a point in $Q_T$. %(convergence on this kind of subsequence was first discussed in \cite{hannah2016unbounded}).
%
%%%%%%%%%%%%%%%%%%%%%%%%%%%%%%%%%%%%%%%%
\begin{lemma}[sufficient descent for unbounded deterministic delays]\label{lemma:deterministic-unbounded-descent}
%%%%%%%%%%%%%%%%%%%%%%%%%%%%%%%%%%%%%%%%
\textbf{Conditions:} Let $f$ be a function (which may be nonconvex) with $L$-Lipschitz gradient and finite $\min f$. The delays $j(k)$ are deterministic and obey $\liminf{j(k)}<\infty$. Use step size $\gamma_k=c/ D_{j(k)}$ for arbitrary fixed $0<c<1$.
\textbf{Results:}
We have
\small
\begin{eqnarray}\label{deterministic-unbounded-descent-result-1}
H_k-H_{k+1}\geq L(\tfrac{1}{\gamma_k}-D_{j(k)})\|\Delta^k\|_2^2
\end{eqnarray}
\normalsize
and
\small
\begin{align}\label{eq:deterministic-unbounded-Dela-to-0}
  \lim_{k}\|\Delta^k\|_2=0.
\end{align}
\normalsize
On any subsequence $Q_T$ (for arbitrarily large $T$), we have:
\small
\begin{align*}
    \lim_{(k\in Q_T)\to\infty}\|d^k\|_2=0,\quad \lim_{(k\in Q_T)\to\infty}\|\nabla_{i_k} f(\hat{x}^{k})\|_2=0,
\end{align*}
\normalsize
\end{lemma}
%%%%%%%%%%%%%%%%%%%%%%%%%%%%%%%%%%%%%%%%

%%%%%%%%%%%%%%%%%%%%%%%%%%%%%%%%%%%%%%%%
%%%%%%%%%%%%%%%%%%%%%%%%%%%%%%%%%%%%%%%%
To prove our next result, we need a new assumption: essentially cyclically semi-unbounded delay (ECSD), which is slightly stronger than the essentially cyclic assumption: in every window of $N'$ steps, every index $i$ is updated at least once with a delay less than $B$ (at iteration $K(i,t)$). The number $B$ just needs to exist and can be arbitrarily large. It does not affect the step size.

%%%%%%%%%%%%%%%%%%%%%%%%%%%%%%%%%%%%%%%%
\begin{theorem}\label{thm:deterministic-unbounded-deterministic}
%%%%%%%%%%%%%%%%%%%%%%%%%%%%%%%%%%%%%%%%
Let the conditions of Lemma \ref{lemma:deterministic-unbounded-descent} hold. For the deterministic index rule under the ECSD assumption, for $T\geq B$, we have:
\small
\begin{align}\label{deterministic-unbounded-deterministic-result-1}
    \lim_{(k\in Q_T)\to\infty}\|\nabla f(x^k)\|_2=0.
\end{align}
\normalsize

\end{theorem}
%%%%%%%%%%%%%%%%%%%%%%%%%%%%%%%%%%%%%%%%

%%%%%%%%%%%%%%%%%%%%%%%%%%%%%%%%%%%%%%
\section{Conclusion}
%%%%%%%%%%%%%%%%%%%%%%%%%%%%%%%%%%%%%%%%
In summary, we have proven a selection of convergence results for async-BCD under bounded and unbounded delays, and stochastic and deterministic block choices. These results do not require the independence assumption that occurs in the vast majority of other work so far. These results were obtained with the use of Lyapunov function techniques, and treating delays directly, rather than modelling them as noise. Future work may involve obtaining a more exhaustive list of convergence results, sharper convergence rates, and an extension to asynchronous stochastic gradient descent-like algorithms. %, such as SAGA.
%%%%%%%%%%%%%%%%%%%%%%%%%%%%%%%%%%%%%%

\newpage
%\section*{Supplementary materials for \textit{Asynchronous Coordinate Descent under More  Realistic Assumptions}}
\section*{Appendix}
Our analysis uses the following standard inequalities. For any $x^1,x^2,\ldots,x^M\in \textrm{R}^N$ and $\varepsilon>0$, it holds that
\small
\begin{align}
\langle x^1,x^2\rangle&\leq \varepsilon\|x^1\|_2^2+\frac{1}{\varepsilon}\|x^2\|_2^2\label{toollemma1}\\
 \langle x^1,x^2\rangle &\leq \|x^1\|_2\cdot\|x^2\|_2\label{toollemma2}\\
 \big\|\sum_{i=1}^M x^i\big\|_2^2 &\leq M\big(\sum_{i=1}^M\|x^i\|_2^2\big)\label{toollemma3}\\
 \|d^k\|_2 &\leq \sum^{k-1}_{i=k-j(k)}\|\Delta^i\|_2\label{toollemma4}
\end{align}
\normalsize
The last inequality is derived from \eqref{eq:Dddef}, where $d^k$ is defined, using a telescoping sum and the triangle inequality. %\commrh{Shoudl we prove this last one?}

\subsection*{Proof of Lemma \ref{lemma:bounded-descent}}
%%%%%%%%%%%%%%%%%%%%%%%%%%%%%%%%%%%%%%%%
 Note that $\Delta^k_{i}=\delta(i,i_k)\cdot\Delta^k_{i_k}$, where %$\delta(i,i_k)$ denotes the Kronecker delta:
 $\delta(i,i_k)=\left\{\begin{array}{ll}
                                      0,&~i=i_k \\
                                      1,&~\textrm{else}
                                    \end{array}
 \right.$.
Recalling the algorithm (\ref{eq:ARock-algorithm}), we have:
\begin{align}\label{bounded-descent-temp-1}
   -\langle \Delta^k,\nabla f(\hat{x}^k)\rangle=-\langle \Delta^k_{i_k},\nabla_{i_k} f(\hat{x}^k)\rangle=\tfrac{L}{\gamma}\| \Delta^k\|^2_2.
\end{align}
Since $\nabla f$ is $L$-Lipschitz,
\begin{align}\label{bounded-descent-temp-2}
    f(x^{k+1})\leq f(x^k)+\langle \nabla f(x^k),\Delta^k\rangle+\frac{L}{2}\|\Delta^k\|^2_2.
\end{align}
Hence
\begin{align}\label{bounded-descent-temp-3}
f(x^{k+1})-f(x^k)&\overset{\eqref{bounded-descent-temp-1}\eqref{bounded-descent-temp-2}}{\leq} \langle \nabla f(x^k)-\nabla f(\hat{x}^k),\Delta^k\rangle +(\tfrac{L}{2}-\tfrac{L}{\gamma})\|\Delta^k\|_2^2\nonumber\\
    &\overset{a)}{\leq} L\|d^k\|_2\cdot\|\Delta^k\|_2+(\tfrac{L}{2}-\tfrac{L}{\gamma})\|\Delta^k\|_2^2\nonumber\\
    &\overset{\eqref{toollemma4}}{\leq} L\sum_{d=k-\tau}^{k-1}\|\Delta^d\|_2\cdot\|\Delta^k\|_2+(\tfrac{L}{2}-\tfrac{L}{\gamma})\|\Delta^k\|_2^2\nonumber\\
    &\overset{b)}{\leq}\tfrac{L}{2\varepsilon}\sum_{i=k-\tau}^{k-1}\|\Delta^i\|_2^2+\left[\tfrac{(\tau\varepsilon+1) L}{2}-\tfrac{L}{\gamma}\right]\|\Delta^k\|_2^2,
\end{align}
where a) follows from \eqref{toollemma2} and the Lipschitz of $\nabla f$, and c) is obtained by applying $a\cdot b\le \tfrac{1}{2\varepsilon} |a|^2+\tfrac{1}{2\varepsilon}|b|^2$ to each term in the sum.
%
%where $\varepsilon>0$ is arbitrary. The first inequality is the directly from (\ref{bounded-descent-temp-1}) and (\ref{bounded-descent-temp-2}). The second is from the Lipschitz continuity of $\nabla f$.  The third one is from the (\ref{toollemma4}).

If $\gamma<\frac{2}{2\tau+1}$, we can choose $\varepsilon>0$ such that $\varepsilon+\frac{1}{\varepsilon}=1+\frac{1}{\tau}(\frac{1}{\gamma}-\frac{1}{2})$.
Then, it can be verified by direct calculation and substitutions that we have:
\begin{align}
    &\xi_k-\xi_{k+1}\overset{\eqref{Lyapunov1}}{=} f(x^k)-f(x^{k+1})+\tfrac{L}{2\varepsilon}\sum_{i=k-\tau}^{k-1}(i-(k-\tau)+1)\|\Delta^i\|_2^2\nonumber\\
    &\quad-\tfrac{L}{2\varepsilon}\sum_{i=k+1-\tau}^{k-1}(i-(k-\tau))\|\Delta^i\|_2^2-\tfrac{L}{2\varepsilon}\tau\|\Delta^k\|_2^2\nonumber\\
    &\quad\overset{c)}{=}f(x^k)-f(x^{k+1})+\tfrac{L}{2\varepsilon}\sum^{k-1}_{i=k-\tau}\|\Delta^i\|_2^2-\tfrac{L}{2\varepsilon}\tau\|\Delta^k\|_2^2\overset{\eqref{bounded-descent-temp-3}}{\geq}\tfrac{1}{2}(\tfrac{1}{\gamma}-\tfrac{1}{2}-\tau)L\cdot\|\Delta^k\|_2^2,
\end{align}
where c) follows from $(i-(k-\tau)+1)\|\Delta^i\|_2^2 - (i-(k-\tau))\|\Delta^i\|_2^2 = \|\Delta^i\|_2^2$.
%the last inequality depends on (\ref{bounded-descent-temp-3}).
%
Therefore we have $\|\Delta^k\|_2^2\in\ell^1$ by using a telescoping sum\footnote{We say a sequence $a^k$ is in $\ell^1$ if $\sum_{k=1}^\infty |a^k|<\infty$.}. This immediately implies (\ref{bounded-descent-result-2}), and (\ref{bounded-descent-result-3}) follows from [Lemma 3, \cite{davis2016convergence}].
\subsection*{Proof of Theorem \ref{thm:bounded-deterministic}}
%%%%%%%%%%%%%%%%%%%%%%%%%%%%%%%%%%%%%%%%
Let $t=t(k)=\lfloor k/N'\rfloor$. Recall $K(i,t)$ is defined at Sec. 1.1. Notice we have:
\begin{align}\label{bounded-deterministic-temp-1}
\|\nabla_i f(x^k)\|_2 &\overset{a)}{\leq} \|\nabla_i f(\hat{x}^{K(i,t)})\|_2 + \|\nabla_i f(x^k)-\nabla_i f(\hat{x}^{k})\|_2\nonumber+\|\nabla_i f(\hat{x}^{k})-\nabla_i f(\hat{x}^{K(i,t)})\|_2\nonumber\\
&\overset{b)}{\leq} \|\nabla_i f(\hat{x}^{K(i,t)})\|_2+L\|d^k\|_2+L\sum_{j=K(i,t)}^{k-1}\|\hat{x}^{j+1}-\hat{x}^j\|_2,
\end{align}
where a) is by the triangle inequality and b) by Lipschitz of $\nabla f$ and then applying the triangle inequality to the expansion of $\|\hat{x}^{k}-\hat{x}^{K(i,t)}\|$.
We now bound each of the right-hand terms.

From  Lemma \ref{lemma:bounded-descent} and by (\ref{toollemma4}), we have
\begin{align}
    \lim_k\|d^k\|_2\leq \lim_k\sum_{i=k-\tau}^{k-1}\|\Delta^i\|_2=0.
\end{align}
By the triangle inequality, we can derive
\begin{align}\label{eq:dklew}
    \|\hat{x}^{k+1}-\hat{x}^{k}\|_2&\leq\|d^k\|_2+\|d^{k+1}\|_2+\|\Delta ^{k}\|_2.
\end{align}
Taking the limitation,
\begin{align}
        \lim_k\|\hat{x}^{k+1}-\hat{x}^{k}\|_2 = 0.
\end{align}
Now notice:
\begin{align}\label{eq:nabf}
\|\nabla_i f(\hat{x}^{K(i,t)})\|_2 = \|\nabla_{i_{K(i,t)}} f(\hat{x}^{K(i,t)})\|_2=\frac{L}{\gamma}\|d^{K(i,t)}\|_2.
\end{align}
Since, as $k\to\infty$, $K(i,t)\to\infty$ and $\|d^k\|_2\to 0$, this last term converges to $0$ and the limit result is proven. The running best rate is obtained through the following argument: since $\|\Delta^k\|_2$ is square summable (by Lemma \ref{lemma:bounded-descent}), so are $\|d^k\|_2$ by (\ref{toollemma4}),  $\|\hat{x}^{k+1}-\hat{x}^{k}\|_2$ by \eqref{eq:dklew}, and $\|\nabla_i f(\hat{x}^{K(i,t)})\|_2$ (in $t=\Theta(k)$) by \eqref{eq:nabf}. Hence, $\|\nabla_i f(x^k)\|_2$ is square summable. With $\lim_k \|\nabla f(x^k)\|_2= 0$, we obtain the running best rate again from [Lemma 3, \cite{davis2016convergence}].
%
%%%%%%%%%%%%%%%%%%%%%%%%%%%%%%%%%%%%%%%%%%%%%%%%%%%%%%%%%%%%%%%%%%%%%%%%%%%%%%%%%%%%
\subsection*{Proof of Theorem \ref{thm:bounded-stochastic-general}}
%%%%%%%%%%%%%%%%%%%%%%%%%%%%%%%%%%%%%%%%
% Under the stochastic rule, we have:\commwy{need to explain this}
% %
% \begin{align}
%  \EE(\|\nabla_{i_k} f(x^{k-\tau})\|_2\mid\chi^{k-\tau})=\frac{\sum_{i=1}^N\|\nabla_{i} f(x^{k-\tau})\|_2}{N}.
% \end{align}
%
Taking the expectation on both sides of \eqref{eq:condEcond} and multiplying $N$ yield
\begin{align}\label{NEnablaik}
 N\EE\|\nabla_{i_k}f(x^{k-\tau})\|_2=\sum_{i=1}^N\EE\|\nabla_{i} f(x^{k-\tau})\|_2.
\end{align}
By $\|\cdot\|_2 \leq \|\cdot\|_1$, we get:
\begin{align}\label{bounded-stochastic-general-temp-7}
\EE \|\nabla f(x^{k-\tau}) \|_2 \leq  \sum_{i=1}^N\EE\|\nabla_i f(x^{k-\tau}) \|_2\overset{\eqref{NEnablaik}}{=} N\EE\|\nabla_{i_k} f(x^{k-\tau})\|_2.
\end{align}
In the next part, we prove $\EE\|\nabla_{i_k} f(x^{k-\tau})\|_2\to 0$.
From (\ref{bounded-descent-result-1}), we can see that $(\|\Delta^k\|_2)_{k\geq 0}$ is bounded. The dominated convergence theorem implies:
\begin{align}\label{bounded-stochastic-general-temp-1}
\lim_k\EE\|\Delta^k\|_2 &= 0.
\end{align}
By (\ref{toollemma4}), we have:
\begin{align}\label{bounded-stochastic-general-temp-3}
    \lim_k \EE( \|d^k\|_2)=0.
\end{align}
Hence,
\begin{align}
\lim_k \EE\|\nabla_{i_k} f(\hat{x}^{k})\|_2\overset{\eqref{eq:ARock-algorithm}}{=}\frac{L}{\gamma}\lim_k \EE\|\Delta^k\|_2=0.\label{bounded-stochastic-general-temp-5}
\end{align}
The triangle inequality and $L$-Lipschitz continuity yield
\begin{align}\label{bounded-stochastic-general-temp-6}
   \EE\|\nabla_{i_k} f(x^{k-\tau})\|_2&\leq\EE\|\nabla_{i_k} f(\hat{x}^k)\|_2+\EE\|\nabla_{i_k} f(x^k)-\nabla_{i_k} f(\hat{x}^k)\|_2\nonumber\\
   &+\EE\|\nabla_{i_k} f(x^k)-\nabla_{i_k} f(x^{k-\tau})\|_2\nonumber\\
&\leq\EE\|\nabla_{i_k} f(\hat{x}^k)\|_2+L\cdot\EE\|d^k\|_2+L\sum_{i=k-\tau}^{k-1}\EE\|\Delta^i\|_2.
\end{align}
Applying \eqref{bounded-stochastic-general-temp-1}, (\ref{bounded-stochastic-general-temp-3}), and (\ref{bounded-stochastic-general-temp-5}) to (\ref{bounded-stochastic-general-temp-6}) yields
\begin{align}\label{bounded-stochastic-general-temp-8}
   \lim_k\EE\|\nabla_{i_k} f(x^{k-\tau})\|_2=0.
\end{align}
With (\ref{bounded-stochastic-general-temp-7}), (\ref{bounded-stochastic-general-temp-8}) yields
\begin{align}\label{bounded-stochastic-general-temp-9}
   \lim_k\EE\|\nabla f(x^{k-\tau})\|_2=0,
\end{align}
which is equivalent to
\begin{align}\label{bounded-stochastic-general-temp-10}
   \lim_k\EE\|\nabla f(x^{k})\|_2=0.
\end{align}
Following a proof similar to that of Theorem \ref{thm:bounded-deterministic}, $\EE\|\nabla f(x^k)\|_2^2$ is summable and thus has the running best rate.

\subsection*{Proof of Lemma \ref{lem:bounded-stochastic-convex}}
%%%%%%%%%%%%%%%%%%%%%%%%%%%%%%%%%%%%%%%%
The proof consists of  two steps: in the first one, we prove
\begin{align}\label{bounded-stochastic-convex-temp-0}
    \pi_k-\pi_{k+1}\geq\frac{L}{4\tau}(\frac{1}{\gamma}-\frac{1}{2}-\tau)\cdot(\EE S(k+1,\tau+1)),
\end{align}
while in the second one, we prove
\begin{align}\label{bounded-stochastic-convex-temp-N2}
   \pi_k^2\leq\beta\cdot (\EE S(k+1,\tau+1))\cdot(\EE S(k,\tau)+\EE\|x^k-\overline{x^k}\|_2^2).
\end{align}
Combining (\ref{bounded-stochastic-convex-temp-0}) and (\ref{bounded-stochastic-convex-temp-N2}) gives us the claim in the lemma.

\textbf{Proving (\ref{bounded-stochastic-convex-temp-0}):} Since $\gamma<\frac{2}{2\tau+1}$, we can choose $\varepsilon>0$ such that
\begin{align}\label{bounded-stochastic-convex-temp-N1}
   \varepsilon+\frac{1}{\varepsilon}=1+\frac{1}{\tau}(\frac{1}{\gamma}-\frac{1}{2})
\end{align}
Direct subtraction of $F_k$ and $F_{k+1}$ yields:
\begin{align}\label{bounded-stochastic-convex-temp-1}
    &F_k-F_{k+1}\overset{a)}{\geq} f(x^k)-f(x^{k+1})+\delta\sum_{i=k-\tau}^{k-1}(i-(k-\tau)+1)\|\Delta^i\|_2^2\nonumber\\
    &\quad-\delta\sum_{i=k+1-\tau}^{k-1}(i-(k-\tau))\|\Delta^i\|_2^2-\delta\tau\|\Delta^k\|_2\nonumber\\
    &\quad\overset{b)}{=} f(x^k)-f(x^{k+1})+\delta S(k,\tau)-\delta\tau\|\Delta^k\|_2 \nonumber\\
    &\quad\overset{c)}{\geq} (\delta-\tfrac{L}{2\varepsilon})S(k,\tau)+\left[\tfrac{L}{\gamma}-\tfrac{(\tau\varepsilon+1) L}{2}-\delta\tau\right]\|\Delta^k\|_2^2\nonumber\\
    &\quad\overset{d)}{=}\tfrac{L}{4\tau}(\tfrac{1}{\gamma}-\tfrac{1}{2}-\tau)\cdot S(k,\tau)+\tfrac{L}{4}(\tfrac{1}{\gamma}-\tfrac{1}{2}-\tau)\cdot\|\Delta^k\|_2^2\nonumber\\
    &\quad\overset{e)}{\geq}\tfrac{L}{4\tau}(\tfrac{1}{\gamma}-\tfrac{1}{2}-\tau)\cdot S(k,\tau)+\tfrac{L}{4\tau}(\tfrac{1}{\gamma}-\tfrac{1}{2}-\tau)\cdot\|\Delta^k\|_2^2\nonumber\\
    &\quad\overset{f)}{=}\tfrac{L}{4\tau}(\tfrac{1}{\gamma}-\tfrac{1}{2}-\tau)\cdot S(k+1,\tau+1),
\end{align}
where a) follows from the definition $F_k$, b) from the definition of $S(k,\tau)$, c) from (\ref{bounded-descent-temp-3}), d) is a direct computation using (\ref{bounded-stochastic-convex-temp-N1}), e) is due to $\tau\geq 1$, and f) is also a result of the definition of $S(k,\tau)$.%\commrh{Need to check last line. But essentially correct.}

\textbf{Proving (\ref{bounded-stochastic-convex-temp-N2}):}
The convexity of $f$ yields
\begin{align}\label{bounded-stochastic-convex-temp-2}
    f(x^k)- f(\overline{x^k})\leq\langle\nabla f(x^k),\overline{x^k}-x^k\rangle.
\end{align}
Let
 \begin{align}
    a^k:=\left(
                              \begin{array}{c}
                                \overline{x^k}-x^k \\
                                \sqrt{\delta\tau}\Delta^{k-1} \\
                               \vdots \\
                                 \sqrt{\delta\tau}\Delta^{k-\tau}\\
                              \end{array}
                            \right),\quad b^k:=\left(
                              \begin{array}{c}
                                \nabla f(x^k) \\
                                \sqrt{\delta\tau}\Delta^{k-1} \\
                               \vdots \\
                                 \sqrt{\delta\tau}\Delta^{k-\tau}\\
                              \end{array}
                            \right).
\end{align}
Using this and the definition of $F_k$ (\ref{Lyapunov2}), we have:
\begin{align}\label{bounded-stochastic-convex-temp-3}
    F_k-\min f\leq \langle a^k,b^k\rangle \leq\|a^k\|_2\|b^k\|_2.
\end{align}
We bound $\EE\|\nabla_{i_k} f(x^{k-\tau})\|_2^2$ as follows:
\begin{align}\label{bounded-stochastic-convex-temp-6}
    \EE\|\nabla_{i_k} f(x^{k-\tau})\|_2^2&
    \overset{a)}{\leq}\EE\big(\|\nabla_{i_k} f(x^{k})\|_2+\|\nabla_{i_k} f(x^{k-\tau})-\nabla_{i_k} f(x^{k})\|_2\big)^2\nonumber\\
    &\overset{b)}{\leq}2\EE\|\nabla_{i_k} f(x^k)\|_2^2+2L^2\tau\sum_{i=k-\tau}^{k-1}\EE\|\Delta^i\|_2^2\nonumber\\
    &\overset{c)}{\leq}4\EE\|\nabla_{i_k} f(\hat{x}^k)\|_2^2+4L^2\EE\|d^k\|_2^2+2L^2\tau\sum_{i=k-\tau}^{k-1}\EE\|\Delta^i\|_2^2\nonumber\\
    &=\tfrac{4L^2}{\gamma^2}\EE\|\Delta^k\|_2^2+6L^2\tau\sum_{i=k-\tau}^{k-1}\EE\|\Delta^i\|_2^2,
\end{align}
where a) follows from the triangle inequality, b) from the Lipschitz of $\nabla f$ and (\ref{toollemma3}), and c) from $\|\nabla_{i_k} f(x^k)\|_2^2\leq 2\|\nabla_{i_k} f(\hat{x}^k)\|_2^2+2\|d^k\|_2^2$ and (\ref{toollemma4}).
We also have the bound
\begin{align}\label{bounded-stochastic-convex-temp+1}
\|\nabla f(x^k)\|_2^2\leq2\|\nabla f(x^{k-\tau})\|_2^2+2L^2\tau\sum_{i=k-\tau}^{k-1}\|\Delta^i\|_2^2,
\end{align}
%By \rev{...(...)...}, the left-hand sides of \eqref{bounded-stochastic-convex-temp-6} and \eqref{bounded-stochastic-convex-temp+1} obey
%\begin{align}
%  \EE(\|\nabla_{i_k}f(x^{k-\tau})\|_2^2\mid\chi^{k-\tau})&=\tfrac{1}{N}\|\nabla f(x^{k-\tau})\|_2^2,\nonumber\\
%   \EE\|\nabla_{i_k} f(x^{k-\tau})\|_2^2&=\tfrac{1}{N}\EE\|\nabla f(x^{k-\tau})\|_2^2\label{bounded-stochastic-convex-temp-5}.
Hence, applying (\ref{NEnablaik}) to (\ref{bounded-stochastic-convex-temp-6}) yields
\begin{align}
    \EE\|\nabla f(x^{k-\tau})\|_2^2\leq\tfrac{4NL^2}{\gamma^2}\EE\|\Delta^k\|_2^2+6NL^2\tau\sum_{i=k-\tau}^{k-1}\EE\|\Delta^i\|_2^2,\nonumber
\end{align}
and further with (\ref{bounded-stochastic-convex-temp+1}),
\begin{align}\label{eenbf22it}
    \EE\|\nabla f(x^{k})\|_2^2\leq\tfrac{8NL^2}{\gamma^2}\EE\|\Delta^k\|_2^2+(12N+2)L^2\tau\sum_{i=k-\tau}^{k-1}\EE\|\Delta^i\|_2^2.
\end{align}
Finally we obtain (\ref{bounded-stochastic-convex-temp-N2}) from
\begin{align}\label{bounded-stochastic-convex-temp-4}
    \pi_k^2=&[\EE(F_k- \min f)]^2\overset{\eqref{bounded-stochastic-convex-temp-3}}{\leq}\EE(\|a^k\|_2\|b^k\|_2)^2\leq\EE(\|a^k\|_2^2)\cdot\EE(\|b^k\|_2^2)\nonumber\\
    &\overset{a)}{\leq} (\tau\EE S(k,\tau)+\EE \|\nabla f(x^k)\|_2^2)\times(\tau\EE S(k,\tau)+\EE\|x^k-\overline{x^k}\|_2^2)\nonumber\\
    &\overset{b)}{\leq} \beta\EE S(k+1,\tau+1)\cdot(\tau\EE S(k,\tau)+\EE\|x^k-\overline{x^k}\|_2^2),
\end{align}
where a) follows from the definitions of $a^k,b^k$ and b) from \eqref{eenbf22it} and the definition of $S(k,\tau)$.
\subsection*{Proof of Theorem \ref{thm:bounded-stochastic-convex}}
%%%%%%%%%%%%%%%%%%%%%%%%%%%%%%%%%%%%%%%%
With (\ref{bounded-stochastic-convex-temp-1}), we can see that $f(x^k)\leq F_k\leq F_0$. Since $f$ is coercive, the sequence $(x^k)_{k\geq 0}$ is bounded. Hence, we have $\sup_{k}\{\|x^k-\overline{x^k}\|_2\}<+\infty$.
Hence, there exists $R>0$ such that
\begin{align}
    \alpha(\sum_{i=k-\tau}^{k-1}\tau\delta\EE\|\Delta^{i}\|_2^2+\EE\|x^k-\overline{x^k}\|_2^2)\leq\tfrac{1}{R}.
\end{align}
for all $k$. Using Lemma \ref{lem:bounded-stochastic-convex}, we have
 \begin{align}
   \pi_k-\pi_{k+1}\geq R \pi_{k}^2.
\end{align}
Using (\ref{bounded-stochastic-convex-temp-1}), we can see that $\pi_k\geq\pi_{k+1}$ for all $k$.
Thus, we have
 \begin{align}
   \pi_k-\pi_{k+1}&\geq R \pi_{k+1}\pi_k\\
   \implies \tfrac{1}{\pi_{k+1}}-\tfrac{1}{\pi_k}&\geq R.
\end{align}
Therefore, using a telescoping sum, we can deduce that:
\begin{align}
    \pi_{k+1}\leq \tfrac{1}{k R+\tfrac{1}{\pi_{0}}}.
\end{align}
Noting $\EE(f(x^k)-\min f)\leq\pi_k$, we have proven the result.
\subsection*{Proof of Theorem \ref{thm:bounded-stochastic-restrictedstronglyconvex}}
%%%%%%%%%%%%%%%%%%%%%%%%%%%%%%%%%%%%%%%%
We have
\begin{align}
    \EE(f(x^k)-\min f)\geq \nu\EE\|x^k-\overline{x^k}\|_2^2,
\end{align}
Hence recalling the definition from (\ref{Lyapunov2}), we have
\begin{align}
    \EE\pi_k\geq \nu\EE\|x^k-\overline{x^k}\|_2^2+\sum_{i=k-\tau}^{k-1}\delta\EE\|\Delta^i\|_2^2\geq\min\{\nu,1\}(\EE\|x^k-\overline{x^k}\|_2^2+S(k,\tau)).\nonumber
\end{align}
Using this, the monotonicity of $\pi^k$, and Lemma \ref{lem:bounded-stochastic-convex} yields
\begin{align}
  \pi_{k}\pi_{k+1}\leq (\pi_{k})^2\leq\tfrac{\alpha }{\min\{\nu,1\}}(\pi_{k}-\pi_{k+1})\cdot \pi_{k}.
\end{align}
Rearranging this yields the result.
%

%%%%%%%%%%%%%%%%%%%%%%%%%%%%%%%%%%%%%%
\subsection*{Proof of Lemma \ref{lemma:stochastic-unbounded-descent}}
%%%%%%%%%%%%%%%%%%%%%%%%%%%%%%%%%%%%%%%%
The Lipschitz continuity of $\nabla f$ yields
\begin{align}
   &f(x^{k+1})-f(x^k)\leq \langle \nabla f(x^k),\Delta^k\rangle+\tfrac{L}{2}\|\Delta^k\|_2^2\nonumber\\
   &\quad \overset{a)}{=}\langle \nabla f(x^k)-\nabla f(\hat{x}^k),\Delta^k\rangle+(\tfrac{L}{2}-\tfrac{L}{\gamma})\|\Delta^k\|_2^2\nonumber\\
   &\quad\leq L\|d^k\|_2\cdot\|\Delta^k\|_2+(\tfrac{L}{2}-\tfrac{L}{\gamma})\|\Delta^k\|_2^2,
\end{align}
where a) is from $-\tfrac{L}{\gamma}\|\Delta^k\|_2^2=\langle \nabla f(\hat{x}^k), \Delta^k\rangle$.
We bound the expectation of $\|d^k\|_2^2$ over the delay and using (\ref{toollemma3}), we have:
\begin{align}\label{stochastic-unbounded-descent-temp-3}
&\EE_{\vec{j}(k)}\big(\|d^k\|_2^2\mid\chi^k\big)\leq\EE_{\vec{j}(k)}\big(\sum_{l=1}^{j(k)}j(k)\|\Delta^{k-l}\|_2^2\mid\chi^k\big)\nonumber\\
&\quad\leq\sum_{j=1}^{+\infty}j p_j\sum_{l=1}^{j}\|\Delta^{k-l}\|_2^2\overset{b)}{=}\sum_{l=1}^{+\infty}(\sum_{j=l}^{+\infty}j p_j)\|\Delta^{k-l}\|_2^2\overset{c)}{\leq} \sum_{i=0}^{k-1}c_{k-i}\|\Delta^i\|_2^2,
\end{align}
where in b), we switched the order of summation in the double sum, and c) uses $\sum_{j=l}^{+\infty}j p_j\leq c_l$.
Taking total expectation $\EE(\cdot)$ on both sides of (\ref{stochastic-unbounded-descent-temp-3}), we obtain
\begin{align}\label{stochastic-unbounded-descent-temp-2}
   \EE\|d^k\|_2^2\leq\sum_{i=0}^{k-1}c_{k-i}\EE\|\Delta^i\|_2^2\overset{d)}{\leq}\sum_{i=0}^{k-1}c_{k-1-i}\EE\|\Delta^i\|_2^2= R(k-1),
\end{align}
where d) is by the fact $(c_i)_{i\geq 0}$ is descending.
Hence:
\begin{align}\label{stochastic-unbounded-descent-temp-1}
     &\EE[f(x^{k+1})-f(x^k)]\leq L\EE\|d^k\|_2\cdot\|\Delta^k\|_2+(\tfrac{L}{2}-\tfrac{L}{\gamma})\EE\|\Delta^k\|_2^2\nonumber\\
    &\quad\leq\tfrac{L}{2\varepsilon}\EE\|d^{k}\|_2^2+\left[\tfrac{(\varepsilon+1) L}{2}-\tfrac{L}{\gamma}\right]\EE\|\Delta^k\|_2^2\nonumber\\
    &\quad\leq\tfrac{L}{2\varepsilon}\sum_{l=1}^{+\infty}(\sum_{j=l}^{+\infty}j p_j)\EE\|\Delta^{k-l}\|_2^2+\left[\tfrac{(\varepsilon+1) L}{2}-\tfrac{L}{\gamma}\right]\EE\|\Delta^k\|_2^2.
\end{align}
Since $\gamma<\frac{2}{2\sqrt{c_0}+1}$, we can choose $\varepsilon>0$ such that
\begin{align}
   \tfrac{1}{2}(\varepsilon+\tfrac{c_0}{\varepsilon})=\tfrac{1}{\gamma}-\tfrac{1}{2}.
\end{align}
With such $\varepsilon$ and (\ref{stochastic-unbounded-descent-temp-1}), direct calculation using the definition of $G^k$ yields (\ref{stochastic-unbounded-descent-result-1}).
When $\gamma<\tfrac{2}{2\sqrt{c_0}+1}$, $\tfrac{L}{2}(\tfrac{1}{\gamma}-\tfrac{1}{2}-\sqrt{c_0})>0$.
From  (\ref{stochastic-unbounded-descent-result-1}), we can see $(R(k))_{k\geq 0}$ is summable (telescoping sum). Thus, we have $\lim_{k}R(k)=0$. Then note  (\ref{stochastic-unbounded-descent-temp-2}) and
\begin{align}
c_0\EE\|\Delta^k\|_2^2\leq\sum_{i=0}^{k}c_{k-i}\EE(\|\Delta^i\|_2^2)=R(k).
\end{align}
Hence then have
\begin{align}
    \lim_{k}\EE(\|d^k\|_2^2)=0,~~\lim_{k}\EE(\|\Delta^k\|_2^2)=0.
\end{align}
%
%%%%%%%%%%%%%%%%%%%%%%%%%%%%%%%%%%%%%%%%

\subsection*{Proof of Theorem \ref{thm:stochastic-unbounded-deterministic}}
%%%%%%%%%%%%%%%%%%%%%%%%%%%%%%%%%%%%%%%%
%Lemma \ref{lemma:stochastic-unbounded-descent} yields $\EE\|d^k\|_2\to 0$ and $\EE\|\Delta^k\|_2\to 0$.
Let $t=t(k)=\lfloor k/N'\rfloor$. Recalling $K(i,t)$ is defined at Sec. 1.1, we have:
\begin{align}\label{add-t1}
\|\nabla_i f(x^k)\|_2 &\overset{a)}{\leq} \|\nabla_i f(\hat{x}^{K(i,t)})\|_2 +\|\nabla_i f(x^{K(i,t)})-\nabla_i f(\hat{x}^{K(i,t)})\|_2 + \|\nabla_i f(x^k)-\nabla_i f(x^{K(i,t)})\|_2\nonumber\\
&\overset{b)}{\leq} \|\nabla_i f(\hat{x}^{K(i,t)})\|_2+L\|d^{K(i,t)}\|_2+L\sum_{j=K(i,t)}^{k-1}\|\Delta^j\|_2,
\end{align}
where a) is by the triangle inequality and b) by the Lipschitz of $\nabla f$ and then applying the triangle inequality to the expansion of $\|x^{k}-x^{K(i,t)}\|$.
We now bound each of the right-hand terms.

Since, as $k\to\infty$, $K(i,t)\to\infty$. With  the Cauchy-Schwarz inequality and (\ref{stochastic-unbounded-descent-result-2}), we have
\begin{align}\label{add-t3}
    \lim_k\EE\|d^{K(i,t)}\|_2\leq \lim_k(\EE\|d^{K(i,t)}\|_2^2)^{\frac{1}{2}}=0.
\end{align}
By $\lim_j\EE\|\Delta^{j}\|_2\leq \lim_j(\EE\|\Delta^{j}\|_2^2)^{\frac{1}{2}}=0$,
\begin{align}\label{add-t2}
    \lim_{k}L\sum_{j=K(i,t)}^{k-1}\EE\|\Delta^j\|_2=0.
\end{align}
Now notice:
\begin{align}
\|\nabla_i f(\hat{x}^{K(i,t)})\|_2 = \|\nabla_{i_{K(i,t)}} f(\hat{x}^{K(i,t)})\|_2=\frac{L}{\gamma}\|d^{K(i,t)}\|_2.
\end{align}
Since $\EE\|d^{K(i,t)}\|_2\to 0$ as $K(i,t)\to\infty$, %this last term converges to $0$, i.e.,
we have
\begin{align}\label{add-t4}
\lim_{k}\EE\|\nabla_i f(\hat{x}^{K(i,t)})\|_2 = 0.
\end{align}
Taking  expectations on both sides of  (\ref{add-t1}), and using (\ref{add-t3}), (\ref{add-t2}) and (\ref{add-t4}), we then prove the result.

%\commwy{Tao: you're lazy here. Please specify similar to what!}
%\commwy{this is not correct because we did not assume $\mathbb{P}(j(k)=0)>0$, so some careful argument is needed.}

\subsection*{Proof of Theorem \ref{thm:stochastic-unbounded-stochastic}}
%%%%%%%%%%%%%%%%%%%%%%%%%%%%%%%%%%%%%%%%
Recall $j(k)$ defined near \eqref{def:jk}. Similar to the bound of $\|d^k\|_2^2$ in \eqref{toollemma4}, we have
\begin{align}
\EE_{\vec{j}(k)}\big(\|x^k-x^{k-j(k)}\|_2^2\mid\chi^{k}\big)\leq\sum_{i=0}^{k-1}s_{k-1-i}\|\Delta^i\|_2^2.
\end{align}
Taking total expectations of both sides yields
\begin{align}\label{stochastic-unbounded-stochastic-temp-0}
\EE\|x^k-x^{k-j(k)}\|_2^2\leq \sum_{i=0}^{k-1}s_{k-1-i}\EE\|\Delta^i\|_2^2.
\end{align}
We have
\begin{align}\label{stochastic-unbounded-stochastic-temp-1}
    \EE\|\nabla_{i_k} f(x^{k-j(k)})\|_2^2&\overset{a)}{\leq}\EE(\|\nabla_{i_k} f(x^{k})\|_2+\|\nabla_{i_k} f(x^{k-j(k)})-\nabla_{i_k} f(x^{k})\|_2)^2\nonumber\\
    &\overset{b)}{\leq}2\EE\|\nabla_{i_k} f(x^k)\|_2^2+2L^2\EE\|x^k-x^{k-j(k)}\|_2^2\nonumber\\
    &\overset{c)}{\leq}4\EE\|\nabla_{i_k} f(\hat{x}^k)\|_2^2+4L^2\EE\|d^k\|_2^2+2L^2\EE\|x^k-x^{k-j(k)}\|_2^2\nonumber\\
    &\overset{d)}{\leq}\tfrac{4L^2}{\gamma^2}\EE\|\Delta^k\|_2^2+6L^2 \sum_{i=0}^{k-1}s_{k-1-i}\EE\|\Delta^i\|_2^2,
\end{align}
where a) follows from the triangle inequality, b) from the Lipschitz of $\nabla f$ and (\ref{toollemma3}), and c) from $\|\nabla_{i_k} f(x^k)\|_2^2\leq 2\|\nabla_{i_k} f(\hat{x})\|_2^2+2\|d^k\|_2^2$ and (\ref{toollemma4}), and d) from (\ref{stochastic-unbounded-stochastic-temp-0}).
Taking total expectation of both sides of assumption (\ref{unassump}) yields
\begin{align}\label{assumptemp+}
    \EE\|\nabla_{i_k} f(x^{k-j(k)})\|_2^2=\frac{\EE\|\nabla f(x^{k-j(k)})\|_2^2}{N}.
\end{align}
By the triangle inequality,
\begin{align}\label{stochastic-unbounded-stochastic-temp-2}
\|\nabla f(x^k)\|_2^2\leq2\|\nabla f(x^{k-j(k)})\|_2^2+2L^2\sum_{i=0}^{k-1}s_{k-1-i}\EE\|\Delta^i\|_2^2.
\end{align}
Hence, combining  (\ref{assumptemp+})
and (\ref{stochastic-unbounded-stochastic-temp-2}) produces
\begin{align}
    \EE\|\nabla f(x^{k-j(k)})\|_2^2\leq\tfrac{4NL^2}{\gamma^2}\EE\|\Delta^k\|_2^2+6NL^2\sum_{i=0}^{k-1}s_{k-1-i}\EE\|\Delta^i\|_2^2;\nonumber
\end{align}
which is substituted into (\ref{stochastic-unbounded-stochastic-temp-2}) to yield
\begin{align}\label{stochastic-unbounded-stochastic-temp-3}
    \EE\|\nabla f(x^{k})\|_2^2\leq\tfrac{8NL^2}{\gamma^2}\EE\|\Delta^k\|_2^2+(12N+2)L^2\sum_{i=0}^{k-1}s_{k-1-i}\EE\|\Delta^i\|_2^2.
\end{align}
By $\sum_{i=0}^{k-1}s_{k-1-i}\leq \sum_{i=0}^{k-1}c_{k-1-i}\EE\|\Delta^i\|_2^2= R(k-1)$ and  (\ref{stochastic-unbounded-descent-result-1}),
\begin{align}
    \lim_k \EE\|\nabla f(x^k)\|_2^2=0.
\end{align}
The proof is completed by applying the Cauchy-Schwarz inequality
\begin{align}
   \EE\|\nabla f(x^k)\|_2\leq (\EE\|\nabla f(x^k)\|_2^2)^{\frac{1}{2}}.
\end{align}
%\commwy{Why do we need the Cauchy-Schwarz inequality?}
%%%%%%%%%%%%%%%%%%%%%%%%%%%%%%%%%%%%%%%%
\subsection*{Proof of Lemma \ref{lem:stochastic-unbounded-stochastic-convex}}
%%%%%%%%%%%%%%%%%%%%%%%%%%%%%%%%%%%%%%%%
This proof is very similar to Lemma \ref{lem:bounded-stochastic-convex} except that $R(k)$ plays the role of $S(k,\tau)$. Let
 \begin{align}
    a^k=\left(
                              \begin{array}{c}
                                \overline{x^k}-x^k \\
                                \sqrt{c_0\bar{\delta}}\Delta^{k-1} \\
                               \vdots \\
                                 \sqrt{c_k\bar{\delta}}\Delta^0\\
                              \end{array}
                            \right),    b^k=\left(
                              \begin{array}{c}
                                \nabla f(x^k) \\
                                \sqrt{c_0\bar{\delta}}\Delta^{k-1} \\
                               \vdots \\
                                 \sqrt{c_k\bar{\delta}}\Delta^0\\
                              \end{array}
                            \right).
\end{align}
Thus, we have
\begin{align}\label{stochastic-unbounded-stochastic-convex-temp-2}
    G_k- \min f\leq\langle a^k,b^k\rangle\leq\|a^k\|_2\|b^k\|_2.
\end{align}
By taking expectations, we get
\begin{align}
   \EE(G_k-\min f)&\leq\EE(\|a^k\|_2\|b^k\|_2)%& \nonumber\\
   %&
   \leq[\EE\|a^k\|_2^2\cdot\EE\|b^k\|_2^2]^{1/2}.
\end{align}

By (\ref{stochastic-unbounded-stochastic-temp-3}) and the definitions of $a^k,b^k,R(k)$,
we get
\begin{align}\label{stochastic-unbounded-stochastic-convex-temp-3}
    &[\EE(G_k- \min f)]^2\leq\EE(\|a^k\|_2^2)\cdot\EE(\|b^k\|_2^2)\nonumber\\
    &\leq (\bar{\delta}R(k)+\EE \|\nabla f(x^k)\|_2^2)\times(\bar{\delta}R(k)+\EE\|x^k-\overline{x^k}\|_2^2)\nonumber\\
    &\leq \bar{\beta}R(k)\times(R(k)+\EE\|x^k-\overline{x^k}\|_2^2).
\end{align}
Finally, from the definition of $\overline{\alpha}$ and Lemma \ref{lemma:stochastic-unbounded-descent}, the theorem follows.
\subsection*{Proof of Theorem \ref{lem:stochastic-unbounded-stochastic-convex-linear}}
%%%%%%%%%%%%%%%%%%%%%%%%%%%%%%%%%%%%%%%%
We have
\begin{align}
    \EE(f(x^k)-\min f)\geq \nu\EE\|x^k-\overline{x^k}\|_2^2,
\end{align}
which also means that
\begin{align}
    \EE(\overline{\delta}R(k)+\|x^k-\overline{x^k}\|_2^2)\leq\max\{1,\frac{1}{\nu}\}\phi_k.
\end{align}
Lemma \ref{lem:stochastic-unbounded-stochastic-convex} yields
\begin{align}
   (\phi_{k})^2\leq\bar{\alpha}\max\{1,\frac{1}{\nu}\}(\phi_{k}-\phi_{k+1})\cdot (\phi_{k})
\end{align}
Note that $\phi_k$ is decreasing, we obtain
\begin{align}
   \phi_{k+1}\leq\bar{\alpha}\max\{1,\frac{1}{\nu}\}(\phi_{k}-\phi_{k+1}).
\end{align}
Then, we have the result by rearrangement.
\subsection*{Proof of Lemma \ref{lemma:deterministic-unbounded-descent}}
%%%%%%%%%%%%%%%%%%%%%%%%%%%%%%%%%%%%%%%%
\begin{align}\label{deterministic-unbounded-descent-temp-1}
     % &
     f(x^{k+1})%\nonumber\\
    &\overset{a)}{\leq} f(x^k)+L\|d^k\|_2\cdot\|\Delta^k\|_2+(\tfrac{L}{2}-\tfrac{L}{\gamma_k})\|\Delta^k\|_2^2\nonumber\\
    &\overset{b)}{\leq} f(x^k)+L\sum_{l=1}^{j(k)}\|\Delta^{k-l}\|_2\cdot\|\Delta^k\|_2+(\tfrac{L}{2}-\tfrac{L}{\gamma_k})\|\Delta^k\|_2^2\nonumber\\
    &\overset{c)}{\leq} f(x^k)+L\sum_{l=1}^{j(k)}(\tfrac{\epsilon_l}{2}\|\Delta^{k-l}\|_2^2+\tfrac{1}{2\epsilon_l}\|\Delta^k\|_2^2)+(\tfrac{L}{2}-\tfrac{L}{\gamma_k})\|\Delta^k\|_2^2\nonumber\\
     &= f(x^k)+\tfrac{L}{2}\sum_{l=1}^{j(k)}\epsilon_l\|\Delta^{k-l}\|_2^2+\tfrac{L}{2}\sum_{l=1}^{j(k)}\tfrac{1}{\epsilon_l}\|\Delta^k\|_2+(\tfrac{L}{2}-\tfrac{L}{\gamma_k})\|\Delta^k\|_2^2 \nonumber\\
     &\overset{d)}\leq f(x^k)+\tfrac{L}{2}\sum_{l=1}^{+\infty}\epsilon_l\|\Delta^{k-l}\|_2^2+\tfrac{L}{2}(1+\sum_{l=1}^{j(k)}\tfrac{1}{\epsilon_l}-\tfrac{2}{\gamma_k})\|\Delta^k\|_2^2.
\end{align}
where a) follows from Lipschitz of $\nabla f$ and definitions of $d^k,\Delta^k$, b) from the triangle inequality, c) from \eqref{toollemma2}, and d) from $j(k)<\infty$.
Then, a direct calculation yields (\ref{deterministic-unbounded-descent-result-1}). Hence (\ref{eq:deterministic-unbounded-Dela-to-0}) follows by summability: $\|\Delta^k\|^2_2\in \ell^1$.
 \begin{align}
    \lim_{k\in Q_T}\|d^k\|_2\leq\sum_{l=k-T}^{k-1}\lim_{l}\|\Delta^l\|_2=0\\
L(\frac{1}{\gamma_k}-D_{j(k)})\|\Delta^k\|_2^2=\frac{c(1-c)}{L D_{j(k)}}\|\nabla_{i_k} f(\hat{x}^{k})\|_2^2.
\end{align}
Therefore,
\begin{align}
\frac{1}{D_{T}}\sum_{k\in Q_T}\|\nabla_{i_k} f(\hat{x}^{k})\|_2^2<\sum_k \frac{\|\nabla_{i_k} f(\hat{x}^{k})\|_2^2}{ D_{j(k)}}<+\infty.
\end{align}
%
%%%%%%%%%%%%%%%%%%%%%%%%%%%%%%%%%%%%%%%%
\subsection*{Proof of Theorem \ref{thm:deterministic-unbounded-deterministic}}
%%%%%%%%%%%%%%%%%%%%%%%%%%%%%%%%%%%%%%%%
%%%%%%%%%%%%%%%%%%%%%%%%%%%%%%%%%%%%%%%%
For  any $T$ and $k\in Q_T$, let $t=t(k)=\lfloor k/N'\rfloor$, and by the triangle inequality:
\begin{align}
&\|\nabla_i f(x^k)\| \leq \|\nabla_i f(x^{K(i,t)})-\nabla_i f(x^k)\|_2\nonumber\\
&\qquad+\|\nabla_i f(\hat{x}^{K(i,t)})-\nabla_i f(x^{K(i,t)})\|_2 + \|\nabla_i f(\hat{x}^{K(i,t)})\|_2.
\end{align}
From  Lemma \ref{lemma:deterministic-unbounded-descent}, we have
\begin{align}
    &\lim_k\|\nabla_i f(x^{K(i,t)})-\nabla_i f(x^k)\|_2\leq \lim_k L\sum_{i=k-N'+1}^{k-1}\|\Delta^i\|_2=0.
\end{align}
Noting $K(i,t)\in Q_T$ by the ECSD assumption, we can derive
\begin{align}
    &\lim_{k}\|\nabla_i f(\hat{x}^{K(i,t)})-\nabla_i f(x^{K(i,t)})\|_2\leq\quad \lim_{k}L\|d^{K(i,t)}\|_2=0.
\end{align}
Now notice by Lemma \ref{lemma:deterministic-unbounded-descent}:
\begin{align*}
\lim_{k}\|\nabla_i f(\hat{x}^{K(i,t)})\|_2 =\lim_{K(i,t)} \|\nabla_{i_{K(i,t)}} f(\hat{x}^{K(i,t)})\|_2=0.
\end{align*}
Since $K(i,t)\to\infty$, this right term converges to $0$ and the result is proven.

\begin{thebibliography}{10}

\bibitem{cannelli2016asynchronous}
Loris Cannelli, Francisco Facchinei, Vyacheslav Kungurtsev, and Gesualdo
  Scutari.
\newblock Asynchronous parallel algorithms for nonconvex big-data optimization:
  Model and convergence.
\newblock {\em arXiv preprint arXiv:1607.04818}, 2016.

\bibitem{cannelli2017asynchronous}
Loris Cannelli, Francisco Facchinei, Vyacheslav Kungurtsev, and Gesualdo
  Scutari.
\newblock Asynchronous parallel algorithms for nonconvex big-data optimization.
  {Part II}: Complexity and numerical results.
\newblock {\em arXiv preprint arXiv:1701.04900}, 2017.

\bibitem{ChowWuYin2016_cyclic}
Yat~Tin Chow, Tianyu Wu, and Wotao Yin.
\newblock Cyclic coordinate update algorithms for fixed-point problems:
  Analysis and applications.
\newblock {\em SIAM Journal on Scientific Computing}, accepted, 2017.

\bibitem{davis2016asynchronous}
Damek Davis.
\newblock The asynchronous palm algorithm for nonsmooth nonconvex problems.
\newblock {\em arXiv preprint arXiv:1604.00526}, 2016.

\bibitem{davis2016convergence}
Damek Davis and Wotao Yin.
\newblock Convergence rate analysis of several splitting schemes.
\newblock In {\em Splitting Methods in Communication, Imaging, Science, and
  Engineering}, pages 115--163. Springer, 2016.

\bibitem{de2015taming}
Christopher~M De~Sa, Ce~Zhang, Kunle Olukotun, and Christopher R{\'e}.
\newblock Taming the wild: A unified analysis of hogwild-style algorithms.
\newblock In {\em Advances in neural information processing systems}, pages
  2674--2682, 2015.

\bibitem{friedman2007pathwise}
Jerome Friedman, Trevor Hastie, Holger H{\"o}fling, Robert Tibshirani, et~al.
\newblock Pathwise coordinate optimization.
\newblock {\em The Annals of Applied Statistics}, 1(2):302--332, 2007.

\bibitem{friedman2010regularization}
Jerome Friedman, Trevor Hastie, and Rob Tibshirani.
\newblock Regularization paths for generalized linear models via coordinate
  descent.
\newblock {\em Journal of statistical software}, 33(1):1, 2010.

\bibitem{hannah2016unbounded}
Robert Hannah and Wotao Yin.
\newblock On unbounded delays in asynchronous parallel fixed-point algorithms.
\newblock {\em arXiv preprint arXiv:1609.04746}, 2016.

\bibitem{lai2013augmented}
Ming-Jun Lai and Wotao Yin.
\newblock Augmented $\ell_1$ and nuclear-norm models with a globally linearly
  convergent algorithm.
\newblock {\em SIAM Journal on Imaging Sciences}, 6(2):1059--1091, 2013.

\bibitem{LeblondPedregosaLacoste-Julien2016_asaga}
R{\'e}mi Leblond, Fabian Pedregosa, and Simon Lacoste-Julien.
\newblock Asaga: {{Asynchronous Parallel Saga}}.
\newblock 2016.

\bibitem{LiuWright2015_asynchronous}
J.~Liu and S.~Wright.
\newblock Asynchronous stochastic coordinate descent: Parallelism and
  convergence properties.
\newblock 25(1):351--376, 2015.

\bibitem{LiuWrightReBittorfSridhar2015_asynchronous}
Ji~Liu, Stephen~J. Wright, Christopher R{\'e}, Victor Bittorf, and Srikrishna
  Sridhar.
\newblock An asynchronous parallel stochastic coordinate descent algorithm.
\newblock 16(1):285--322, 2015.

\bibitem{mania2015perturbed}
Horia Mania, Xinghao Pan, Dimitris Papailiopoulos, Benjamin Recht, Kannan
  Ramchandran, and Michael~I Jordan.
\newblock Perturbed iterate analysis for asynchronous stochastic optimization.
\newblock {\em arXiv preprint arXiv:1507.06970}, 2015.

\bibitem{peng2016arock}
Zhimin Peng, Yangyang Xu, Ming Yan, and Wotao Yin.
\newblock Arock: an algorithmic framework for asynchronous parallel coordinate
  updates.
\newblock {\em SIAM Journal on Scientific Computing}, 38(5):A2851--A2879, 2016.

\bibitem{pengxuyanyin17}
Zhimin Peng, Yangyang Xu, Ming Yan, and Wotao Yin.
\newblock On the convergence of asynchronous parallel iteration with arbitrary
  delays.
\newblock {\em arXiv preprint arXiv:1612:04425}, 2016.

\bibitem{RechtReWrightNiu2011_hogwild}
Benjamin Recht, Christopher Re, Stephen Wright, and Feng Niu.
\newblock Hogwild!: A lock-free approach to parallelizing stochastic gradient
  descent.
\newblock In J.~Shawe-Taylor, R.~S. Zemel, P.~L. Bartlett, F.~Pereira, and
  K.~Q. Weinberger, editors, {\em Advances in {{Neural Information Processing
  Systems}} 24}, pages 693--701. {Curran Associates, Inc.}, 2011.

\bibitem{xu2017asynchronous}
Yangyang Xu.
\newblock Asynchronous parallel primal-dual block update methods.
\newblock {\em arXiv preprint arXiv:1705.06391}, 2017.

\end{thebibliography}
\end{document}